\documentclass[letterpaper,oneside]{amsart}
\usepackage[english]{babel}
\usepackage{tikz-cd}
\usepackage{geometry,amssymb,graphicx,enumerate,accents,xcolor}
\usepackage{comment}
\newcommand{\assign}{:=}
\newcommand{\infixand}{\text{ and }}
\newcommand{\infixiff}{\text{ iff }}
\newcommand{\nobracket}{}

\newcommand{\tmop}[1]{\ensuremath{\operatorname{#1}}}

\newenvironment{enumeratealpha}{\begin{enumerate}[a{\textup{)}}] }{\end{enumerate}}
\theoremstyle{plain}
\newtheorem{theorem}{Theorem}[section]
\newtheorem{corollary}[theorem]{Corollary}
\newtheorem{proposition}[theorem]{Proposition}
\newtheorem{lemma}[theorem]{Lemma}

\newtheorem*{theorem*}{Theorem}
\newtheorem*{lemma*}{Lemma}
\newtheorem*{prop*}{Proposition}
\newtheorem*{cor*}{Corollary}
\newtheorem*{conj*}{Conjecture}

\theoremstyle{definition}
\newtheorem*{definition*}{Definition}
\newtheorem{definition}[theorem]{Definition}

\theoremstyle{remark}
\newtheorem*{remark*}{Remark}
\newtheorem{remark}[theorem]{Remark}
\newtheorem*{example*}{Example}
\newtheorem{example}[theorem]{Example}

\numberwithin{theorem}{section}
\numberwithin{equation}{section} 
\numberwithin{figure}{section} 

\newcommand{\catname}[1]{\mathbf{#1}}

\begin{document}

\title[Model Category Structure for Discrete and Continuous Homotopy]{Cofibration and model category structures for Discrete and Continuous
Homotopy}
\thanks{This material is based in part upon work supported by the US National
Science Foundation under Grant No. DMS-1928930 while the author participated
in a program supported by the Mathematical Sciences Research Institute. The
program was held in the summer of 2022 in partnership with the the Universidad
Nacional Aut{\'o}noma de M{\'e}xico. This work was also supported by the
CONACYT Investigadoras y Investigadores por M{\'e}xico Project \#1076 and by
the grant N62909-19-1-2134 from the US Office of Naval Research Global and the
Southern Office of Aerospace Research and Development of the US Air Force
Office of Scientific Research.}

\author{Antonio Rieser}
\address{Centro de Investigaci{\'o}n en Matem{\'a}ticas, Calle Jalisco s/n,
C.P. 36023, Guanajuato, GTO, M{\'e}xico}
\email{antonio.rieser@cimat.mx}

\begin{abstract}
  We show that the categories $\catname{PsTop}$ and $\catname{Lim}$ of
  pseudotopological spaces and limit spaces, respectively, admit cofibration
  category structures, and that $\catname{PsTop}$ admits a model
  category structure, giving several ways to simultaneously study the homotopy
  theory of classical topological spaces, combinatorial spaces such as graphs
  and matroids, and metric spaces endowed with a privileged scale, in addition
  to spaces of maps between them. In the process, we give a sufficient
  condition for a topological construct which contains compactly generated
  Hausdorff spaces as a subcategory to admit an $I$-category structure. We
  further show that, for a topological space $X \in \mathcal{C}$, the homotopy
  groups of $X$ constructed in the cofibration category on
  $\catname{PsTop}$ are isomorphic to those constructed classically in
  $\catname{Top^{\ast}}$.
\end{abstract}

{\maketitle}

\section{Introduction}

In {\cite{Steenrod_1967}}, Steenrod proposed that the category of compactly
generated Hausdorff spaces should be the appropriate setting in which to study
the homotopy theory of topological spaces, arguing that this category
contained the principal spaces of interest in algebraic topology, that it was
closed under a number of standard operations, and that the operations in
question satisfy a number of reasonable properties. In modern terminology,
these last two conditions are encapsulated by the fact that the category of
compactly generated Hausdorff spaces is Cartesian closed, and, over time, the
category of compactly generated Hausdorff spaces did indeed become one of the
standard categories in which to study the algebraic topology of (most)
topological spaces.

In recent years, there has been an increasing interest in using tools from
algebraic topology to study combinatorial objects such as graphs
{\cite{Carranza_Kapulkin_arXiv_2022,Babson_etal_2006,Barcelo_Laubenbacher_2005,Plaut_Wilkins_2013,Conant_etal_2014,Grigoryan_etal_2014,Chih_Scull_2021,Chih_Scull_2022}},
as well as an interest in studying the algebraic topology of metric spaces
which have been decorated with a privileged scale
{\cite{Plaut_Wilkins_2013,Conant_etal_2014,Smale_Smale_2012,Bartholdi_etal_2012,Rieser_2021}},
which we call \emph{scaled metric spaces}, following {\cite{Rieser_2021}}.
Furthermore, in topological data analysis {\cite{Carlsson_2009}}, one
encounters the problem of how to compare the topological invariants of a
topological space with the invariants of an object in one or more of the above
classes, making it important to work in a category which simultaneously
contains both combinatorial objects as well as more classical topological
spaces. We therefore find ourselves once again in search of a ``convenient''
category for the study of the algebraic topology of certain spaces of
interest, this time combinatorial objects and scaled metric spaces, in
addition to CW complexes. In {\cite{Rieser_2021}}, we proposed that this could
be at least partially resolved by developing algebraic topology in the
category $\catname{Cl}$ of Cech closure spaces, since this category
contains graphs, scaled metric spaces, and topological spaces as
subcategories, and it also admits non-trivial maps between them. This
perspective was further developed in
{\cite{Bubenik_Milicevic_2021,Bubenik_Milicevic_2021b}}. Unfortunately,
however, the category of Cech closure spaces has a serious shortcoming: like
the category of topological spaces, it is not Cartesian closed, and, in
particular, there is no canonical closure structure on function spaces for
which the evaluation map is continuous. We are also unaware of a subcategory
of $\catname{Cl}$ which is Cartesian closed and which also contains compactly
generated topological spaces, graphs, and scaled metric spaces as full
subcategories, leaving us to investigate Cartesian closed categories which
contain $\catname{Cl}$ itself as a subcategory.

The goal of the present article is to show that the category
$\catname{PsTop}$ of pseudotopological spaces - the Cartesian closed
hull of $\catname{Cl}$ {\cite{Beattie_Butzmann_2002,Preuss_2002}} - can serve
as a convenient setting for the simultaneous study of homotopy theory on many
types of combinatorial spaces as well as CW-complexes, and, in particular, that
it is well-suited to studying homotopy theory on graphs, scaled metric spaces,
topological spaces, and spaces of maps between them. We do this in a number of
steps. First, we show that any Cartesian closed category  C  which contains
compactly generated Hausdorff spaces as a subcategory and satisfies some
reasonable additional conditions admits a cofibration category structure (in
the sense of Baues {\cite{Baues_1989,Baues_1999}}). Second, we use the
cofibration category structure to construct homotopy groups for
pseudotopological spaces, and we show that the resulting homotopy groups are
isomorphic to the classical homotopy groups when the pseudotopological space
is also topological, i.e. that $\pi_n (X) = [S^n, X]$ as sets, and this can be
given a group structure for $n \geq 1$. While the cofibration category
structure is enough to build a rich homotopy theory on $\catname{PsTop}$, we
further improve on this in Section \ref{sec:Model category} by showing that
the homotopy groups from the cofibration category may be used to define the
build a model category structure on $\catname{PsTop}$ which extends the
Quillen model category structure on $\catname{Top}$.

\section{Pseudotopological Spaces}

In this section, we introduce the category $\catname{PsTop}$ of
pseudotopological spaces and continuous functions, and we recall a number of
basic results about $\catname{PsTop}$ and several related categories for later
reference. We also show that pseudotopological spaces contain as full
subcategories the main categories of interest to combinatorics and applied
topology: graphs, scaled metric spaces, matroids, topological spaces, and Cech
closure spaces.

We begin by defining general convergence spaces, specializing this to
pseudotopological spaces in Definition \ref{def:PsTop}

\begin{definition}
  A \emph{convergence space} is a pair $(X, \Lambda)$, where $X$ is a set
  and $\Lambda \subset \mathcal{\mathbb{F}} (X) \times X$ is a relation
  between the set of filters on $X$, which we denote by $\mathbb{F} (X)$, and
  the elements of $X$, where the relation $\Lambda$ satisfies the following
  axioms.
  \begin{enumerate}
    \item If $(\lambda, x) \in \Lambda$ and $\lambda \subset \lambda'$, then
    $(\lambda', x) \in \Lambda$
    
    \item $\dot{x} \in \Lambda$ for every $x \in X$, where $\dot{x}$ denotes
    the filter generated by $x$.
  \end{enumerate}
  We call such a relation $\Lambda$ a \emph{convergence structure} or
  \emph{convergence} on $X$.
  
  If, in addition,
  \begin{enumerate}
    \item $(\lambda, x), (\lambda', x) \in \Lambda$, then $(\lambda \cap
    \lambda', x) \in \Lambda$,
  \end{enumerate}
  then we say that $(X, \Lambda)$ is a \emph{limit space}.
  
  Given two convergence structures $\Lambda$ and $\Lambda'$ on a set $X$, we
  say that $\Lambda$ is \emph{coarser} than $\Lambda'$, and $\Lambda'$ is
  \emph{finer} than $\Lambda$, iff $\Lambda \subset \Lambda'$.
  
  We will often write $\lambda \rightarrow x$ in $(X, \Lambda)$ for $(\lambda,
  x) \in \Lambda$, and we define
  \[ \lim_{\Lambda} \lambda \assign \{ x \in X \mid (\lambda, x) \in \Lambda
     \} . \]
  When writing $\lim_{\Lambda} \lambda$ and $(X, \Lambda)$, we will sometimes
  omit the convergence structure $\Lambda$ when it is unambiguous. Finally,
  given a convergence structure $\Lambda$, we will abuse notation and write
  $\lambda \in \Lambda$ to indicate that there exists and element $x \in X$
  such that $(\lambda, x) \in \Lambda$.
\end{definition}

\begin{definition}
  Let $(X, \Lambda_X)$ and $(Y, \Lambda_Y)$ be convergence spaces, and let $f
  : X \rightarrow Y$ be a map from $X$ to $Y$. We say that $f$ is
  \emph{continuous} iff
  \[ \lambda \rightarrow x \in \Lambda_X \Rightarrow f (\lambda) \rightarrow f
     (x) \in \Lambda_Y, \]
  where $f (\lambda)$ is the filter in $Y$ generated by the collection $\{ f
  (U) \subset Y \mid U \in \lambda \}$ of subsets of $Y$.
\end{definition}

\begin{definition}
  An \emph{ultrafilter} $\gamma$ on a set $X$ is a filter such that, for
  every subset $A \subset X$, either $A \in \gamma$ or $X \backslash A \in
  \gamma$.
\end{definition}
\begin{definition}
  Given a filter $\lambda$ on a set $X$, we denote by $\beta (\lambda)$ the
  set of all ultrafilters containing $\lambda$.
\end{definition}

With this definition, we recall the following.
\begin{proposition}
  \label{lem:Filter intersection of ultrafilters}Any filter $\lambda$ on $X$
  is equal to the intersection of all of the ultrafilters which contain it,
  i.e. $\lambda = \bigcap_{\gamma \in \beta (\lambda)} \gamma$.
\end{proposition}
In particular, Lemma \ref{lem:Filter intersection of ultrafilters} and the
definition of a convergence space implies that $\lim_{\Lambda} \lambda \subset
\bigcap_{\gamma \in \beta (\lambda)} \lim_{\Lambda} \gamma$ for every filter
$\lambda$ on a convergence space. We define a pseudotopological space to be a
convergence space in which the reverse inclusion also holds.

\begin{definition}
  \label{def:PsTop}We say that a convergence space $(X, \Lambda)$ is a
  \emph{pseudotopological space} iff
  \[ \lim_{\Lambda} \lambda = \bigcap_{\gamma \in \beta (\lambda)}
     \lim_{\Lambda} \gamma, \]
  for every filter $\lambda$ on $X$, where $\beta (\lambda)$ is the collection
  of ultrafilters which contains $\lambda$.
\end{definition}

Finally, we recall from {\cite{Preuss_2002}} that pseudotopological spaces are
also limit spaces.

\begin{proposition}[{\cite{Preuss_2002}}, Remark 2.3.1.2]
  The category $\catname{PsTop}$ of pseudotopological spaces and continuous
  maps between them is a full subcategory of the category $\catname{Lim}$ of
  limit spaces and continuous maps.
\end{proposition}

\subsection{Classes of Pseudotopological Spaces}

In this section, we give the principal examples of pseudotopological spaces
which are of interest to discrete and continuous homotopy.

\subsubsection{Pseudotopologies from Cech Closure Spaces}

We recall from {\cite{Cech_1966}} that a Cech closure space (also called a
\emph{pretopological space}) is a pair $(X, c)$ where $X$ is a set and $c
: \mathcal{P} (X) \rightarrow \mathcal{P} (X)$ is a function on the power set
of $X$ which satisfies
\begin{enumerate}
  \item $c (\emptyset) = \emptyset$,
  
  \item $A \subset c (A)  \forall A \subset X$, and
  
  \item $c (A \cup B) = c (A) \cup c (B)  \forall A, B \subset X$.
\end{enumerate}
A function $f : (X, c_X) \rightarrow (Y, c_Y)$ between closure spaces is said
to be \emph{continuous} iff for all $A \subset X$, $f (c_X (A)) \subset
c_Y (f (A))$. Cech closure spaces (pretopologies) and continuous maps form a
category, denoted $\catname{Cl}$ (or $\catname{PreTop}$).

The fundamental fact about about pretopologies is that they are completely
determined by the choice of a neighborhood filter at every point $x \in X$,
i.e.

\begin{theorem}[{\cite{Cech_1966}}, Theorem 14.B.10]
  Let $X$ be a set, and suppose that, for each element $x \in X$, $\mathcal{U}
  (x)$ is a collection of subsets of $X$ satisfying
  \begin{enumerate}
    \item $\mathcal{U} (x) \neq \emptyset$,
    
    \item For each $U \in \mathcal{U} (x)$, $x \in U$, and
    
    \item For each $U_1$ and $U_2$ in $\mathcal{U} (x)$, there exists a $U \in
    \mathcal{U} (x)$ with $U \subset U_1 \cap U_2$.
  \end{enumerate}
  Then there exists a unique closure operation $c : \mathcal{P} (X)
  \rightarrow \mathcal{P} (X)$ such that, for each $x \in X$, $\mathcal{U}
  (x)$ is a local base for the neighborhood filter of $x \in (X, c)$.
  Furthermore, $c$ is given by the formula
  \[ c (A) : = \{ x \mid x \in X, U \in \mathcal{U} (x) \Rightarrow U \cap X
     \neq \emptyset \} . \]
\end{theorem}

\begin{definition}
  Given a closure space $(X, c)$, we define a convergence structure $\Lambda_c
  \subset \mathbb{F} (X) \times X$ by: $(\lambda, x) \in \Lambda_c$ iff
  $\lambda$ contains the neighborhood filter of $x$.
\end{definition}

\begin{definition}
  Let $(X, c)$ be a closure space. A set $A \subset X$ is said to be
  \emph{closed} iff $c (A) = A$ and $A \subset X$ is said to be
  \emph{open} if $X \backslash A$ is closed.
\end{definition}

\begin{proposition}[{\cite{Preuss_2002}}, Remark 2.3.1.2]
  The maps $(X, c) \mapsto (X, \Lambda_c)$ and $f \mapsto f$, where $f : (X,
  c) \rightarrow (Y, c')$ is a continuous map, defines a functor $F :
  \catname{Cl} \rightarrow \catname{PsTop}$ which makes the
  category $\catname{Cl}$ of Cech closure spaces into a full subcategory of
  $\catname{PsTop}$.
\end{proposition}

Note that, in particular, the above proposition implies that for a closure
structure $c$, the induced convergence structure $\Lambda_c$ is
pseudotopological, and that a continuous function for the closure structure is
continuous for the pseudotopological structure. We now make precise what it
means for a pseudotopological space to be topological.

\begin{definition}
  A pseudotopological space $(X, \Lambda)$ is said to be topological iff there
  exists a topological closure space $(X, c_{\tau})$ for which $\Lambda =
  \Lambda_{\tau}$, i.e. there is a topology $\tau$ on $X$ such that a filter
  $\lambda \rightarrow x$ in $\Lambda$ iff $\lambda$ contains the filter of
  open sets $\mathcal{U} (x)$. 
\end{definition}

Closure structures are extremely common. Every topological space $(X, \tau)$
admits a closure structure $c_{\tau}$, also known as the Kuratowski closure
structure, where $c_{\tau} (A) = \bar{A}$, where $\bar{A}$ denotes the
topological closure of $A$ (see {\cite{Cech_1966}} or {\cite{Rieser_2021}} for
further details). Closure structures on graphs and semi-pseudometric spaces
were discussed in {\cite{Cech_1966}}, Chapter 14, as well as in
{\cite{Rieser_2021}}, and closure structures on scaled metric spaces may be
constructed from an induced semi-pseudometric as described in
{\cite{Rieser_2021}}. Matroids, furthermore, may also be defined using a
(topological) closure structure, which then also induces a pseudotopological
structure.

A closure structure may be induced by a hypergraph in two ways. First, we may
construct a closure structure from the induced graph, where two vertices are
an edge in the graph if they are contained in an edge of the hypergraph. This
closure structure is given by $c (v) : = \{ w \mid \exists e \in E \tmop{such}
\tmop{that} \{ v, w \} \subset e \}$, and $c (A) \assign \cup_{v \in A} c
(v)$, where $v \in V$ and $A \subset V$. We denote the resulting
pseudotopological space by $(V, \Lambda_c)$. Alternately, given a hypergraph
$H = (V, E)$, we may create a topological closure structure on its induced
simplicial complex $\Sigma = \cup_{e \in E} (\mathcal{P} (e) \backslash \{
\emptyset \})$, the union of nonempty subsets of each edge, by taking the
Alexandrov topology on $\Sigma$. Note, however, that neither of these
procedures of constructing closure structures on hypergraphs induces a fully
faithful functor from hypergraphs to closure spaces.

\subsubsection{Pseudotopologies on Spaces of Continuous Functions}

Unlike in $\catname{Top}$ or $\catname{Cl}$, spaces of functions between
pseudotopological spaces have a canonical pseudotopological structure which
makes the evaluation map continuous, called the \emph{continuous
convergence structure.} We recall its definition from
{\cite{Beattie_Butzmann_2002}}.

\begin{definition}
  Given convergence spaces $(X, \Lambda_X)$ and $(Y, \Lambda_Y)$, we denote by
  $\mathcal{C} (X, Y)$ the set of all continuous mappings from $X$ to $Y$, and
  we denote by
  \[ \omega_{X, Y} : \mathcal{C} (X, Y) \times X \rightarrow Y \]
  the evaluation mapping, defined by $\omega_{X, Y} (f, x) \assign f (x)$.
  
  The \emph{continuous convergence structure} $\Lambda_{\mathcal{C}}$ on
  $\mathcal{C} (X, Y)$ is defined by the following:
  \[ \mathcal{H} \rightarrow f \tmop{in} \Lambda_{\mathcal{C}} \infixiff
     \omega_{X, Y} (\mathcal{H} \times \lambda) \rightarrow f (x) \forall
     (\lambda, x) \in \Lambda_X . \]
\end{definition}

We see that, by definition, the evaluation map $\omega_{X, Y}$ is continuous.
Furthermore, when the spaces $X$ and $Y$ are pseudotopological, then
$(\mathcal{C} (X, Y), \Lambda_{\mathcal{C}})$ is also pseudotopological, as
guaranteed by the following proposition.

\begin{proposition}[{\cite{Beattie_Butzmann_2002}}, Theorem 1.5.5]
  $(\mathcal{C} (X, Y), \Lambda_{\mathcal{C}})$ is a pseudotopological space
  iff $(Y, \Lambda_Y)$ is pseudotopological.
\end{proposition}

We note that the continuous convergence structure $\Lambda_{\mathcal{C}}$ is
not in general pretopological
{\cite{Beattie_Butzmann_2002,Dolecki_Mynard_2016}}.
\section{Convenient Categories}

Pseudotopological spaces will be our prototype for a \emph{convenient
category}, i.e. a category which contains the main spaces of interest and
which is Cartesian closed. Before setting the definition of convenience, we
quickly recall several facts about Cartesian closed categories. For more
details, see the relevant chapters in {\cite{Adamek_etal_2006,Riehl_2016}} on
Cartesian closed categories and the books
{\cite{Preuss_2002,Dolecki_Mynard_2016}} for more on the categorical
properties of convergence spaces, limit spaces, and pseudotopological spaces.

\subsection{Cartesian Closed Categories and Convenience}

\begin{definition}
  A category $\mathcal{C}$ is said to be \emph{Cartesian closed} if it has
  all finite products and for each $\mathcal{C}$-object $C$, the functor $(X
  \times -) : \mathcal{C} \rightarrow \mathcal{C}$ is a left adjoint.
\end{definition}

The following well-known proposition is now immediate from the definition.

\begin{proposition}
  Let $\mathcal{C}$ be a Cartesian closed category. Then for all
  $\mathcal{C}$-objects $C$, the functor $(C \times -) : \mathcal{C}
  \rightarrow \mathcal{C}$ preserves colimits.
\end{proposition}

\begin{proof}
  Left adjoints preserve colimits ({\cite{Riehl_2016}}, Theorem 4.5.3).
\end{proof}

\begin{remark}
  Convergence spaces, limit spaces, and pseudotopological spaces are all
  Cartesian closed ({\cite{Preuss_2002}}, 3.1.9(5)(6), 3.3.4(1)). Topological
  spaces and closure spaces, however, are not Cartesian closed
  ({\cite{Preuss_2002}}, 3.1.9(1)).
\end{remark}

For the purposes of the present article, we define a convenient category as
follows.

\begin{definition}
  A \emph{convenient category of spaces} $\mathcal{C}$ is a Cartesian
  closed category which contains the categories of graphs, scaled metric
  spaces, and compactly generated Hausdorff spaces as subcategories.
\end{definition}

\begin{example}
  Generalized convergence spaces, limit spaces, and pseudotopological spaces
  are all convenient categories.
\end{example}

\subsection{Homotopies and the Homotopy Gluing Property}

We now discuss homotopies of maps and homotopy equivalence. In particular, we
present a homotopy gluing criterion which we will see provides a sufficient
condition for a convenient category to
admit an $I$-category structure. We also show that pseudotopological spaces
and limit spaces satisfy this criterion.

\begin{definition}
  \label{def:Cylinder}We say that a quadruple $(I, i_0, i_1, p)$ is a
  \emph{cylinder on a category $\mathcal{C}$} iff $I : \mathcal{C}
  \rightarrow \mathcal{C}$ is a functor and $i_0, i_1 :
  \tmop{Id}_{\mathcal{C}} \rightarrow I$ and $p : I \rightarrow
  \tmop{Id}_{\mathcal{C}}$ are natural transformations such that, for all
  $\mathcal{C}$-objects $X$, the compositions $p \circ i_k : X \rightarrow I X \rightarrow X, k \in \{ 0, 1 \}$
  are the identity. We call the $i_k$ the inclusions of the cylinder and $p$
  the projection of the cylinder. When the $i_k$ and $p$ are understood, we
  will simply refer to the cylinder $(I, i_0, i_1, p)$ by the functor $I$.
\end{definition}
A cylinder on a category $\mathcal{C}$ gives rise to a notion of homotopy in the following way.
\begin{definition}
  \label{def:Homotopy I}Let $f, g : X \rightarrow Y$ be morphisms in a
  category $\mathcal{C}$ with a cylinder $(I, i_0, i_1, p)$. We say that a
  morphism $H : I X \rightarrow Y$ is a \emph{homotopy from $f$ to $g$}
  iff $H \circ i_0 (X) = f$ and $H \circ i_1 (X) = g$.
\end{definition}

As we will see, it may occur that homotopies do not ``glue'' properly, that
is, that homotopy as defined above is not necessarily a transitive relation.
We will require our convenient categories to satisfy one additional
requirement, which, as we will see in Proposition \ref{prop:Homotopy
transitive} below, gives a sufficient condition for the homotopy relation in
Definition \ref{def:Homotopy I} to be transitive. This will be used in the
construction of the $I$-category structure on $\catname{Lim}$ and
$\catname{PsTop}$ in Section \ref{subsec:I-categories}.

\begin{definition*}
  Let $(I, i_0, i_1, p)$ be a cylinder on a category $\mathcal{C}$. Denote by
  $I A \cup_A I A$ the pushout
  \begin{equation*}
  \begin{tikzcd}
     A \ar[r,"i_0"] \ar[d,"i_1"] & IA \ar[d,"q_1"]\\
     IA \ar[r,"q_0"] & IA \cup_A IA
     \end{tikzcd}
  \end{equation*}
  We say that a cylinder $(I, i_0, i_1, p)$ has the \emph{homotopy gluing
  property} iff, for any object $A \in \mathcal{C},$ there exists an
  isomorphism $u : I A \cup_A I A \cong I A$ such that $u \circ q_k \circ i_k
  = i_k$ for $k = 0, 1$, i.e. such that the diagram
  \begin{equation*}
  \begin{tikzcd}
     A \ar[r,"i_k"] \ar[d,"i_k"] & \ar[d,"q_k"] IA \\
     IA & IA \cup_A IA \ar[l,"\cong"',"u"]
  \end{tikzcd}
\end{equation*}
  commutes for $k = 0, 1$.
\end{definition*}

Note that the $i_k$, $k = 0, 1$, in the second diagram above are no longer the
maps which are identified in the pushout, but the other `ends' of the
cylinder.

We now establish the transitivity of the homotopy relation for a cylinder which satisfies the homotopy gluing property.

\begin{proposition}
  \label{prop:Homotopy transitive}Let $f, g, h : X \rightarrow Y$ be morphisms
  in a category  C  with cylinder $(I, i_0, i_1, p)$. Suppose that $F : I X
  \rightarrow Y$ is a homotopy from $f$ to $g$ and that $G : I X \rightarrow
  Y$ is a homotopy from $g$ to $h$. If $(I, i_0, i_1, p)$ has the homotopy
  gluing property, then the map $F \star G \assign \alpha \circ u^{- 1}
  :IX{\rightarrow}Y$ in the commutative diagram below is a homotopy from $f$
  to $h$.
  \begin{equation}
  	\label{diag:Homotopy}
    \begin{tikzcd}    
    X \ar[rrrrdddd,bend right = 100,"g"']
    \ar[rrrrdddd,bend left=50,"g"] \ar[r,"i_0"] \ar[d,"i_1"] & IX \ar[d,"q_1"] \ar[rrddd,bend left=20,"q_1"'] \ar[rrrdddd,bend left=20,"G"] & & & \\
    IX \ar[r,"q_0"]\ar[rrrdd,bend
    right=10,"q_0"]\ar[rrrrddd,bend right=20,"F"'] & IX
    \cup_X IX \ar[dr,"\cong"',"u"] & & & \\    
    & & IX \ar[dr,"u^{-1}","\cong"'] & & \\
    & & & IX \cup_X IX \ar[dr,dashed,"\alpha"] & \\
    & & & & Y
    \end{tikzcd}
  \end{equation}
\end{proposition}
\begin{proof}
  Since $\mathcal{C}$ has the homotopy gluing property with respect to $I$,
  the isomorphisms $I X \cup_X I X \rightarrow I X$ and $I X \rightarrow I X
  \cup_X I X$ exist by hypothesis. Furthermore by Diagram \ref{diag:Homotopy},
  we have
  \begin{eqnarray*}
    F \star G \circ i_0 (X) & = & F \star G \circ u \circ q_0 \circ i_0 (X) =
    \alpha \circ q_0 \circ i_0 = F \circ i_0 = f, \infixand\\
    F \star G \circ i_1 (X) & = & F \star G \circ u \circ q_1 \circ i_1 (X) =
    \alpha \circ q_1 \circ i_1 = G \circ i_1 = h.
  \end{eqnarray*}
  Therefore, $F \star G$ is a homotopy from $f$ to $h$, as desired.
\end{proof}

\subsubsection{Limit Spaces}The standard cylinder $I X \assign X \times [0,
1]$, $I f \assign (f (x), t)$, where $(X, \Lambda)$ is a limit or
pseudotopological space and $([0, 1], \tau)$ is the topological interval with
its topological convergence structure $\tau$, satisfies the homotopy gluing
property in both limit spaces and pseudotopological spaces. Note that $I X
\cup_X I X \cong X \times ([0, 1] \cup_{1 \sim 0} [0, 1])$ since both limit
spaces and pseudotopological spaces are Cartesian closed categories, so to
check that $I X \cup_X I X \cong I X$, for any limit or pseudotopological
space $X$, it suffices to show that $[0, 1] \cup_{1 \sim 0} [0, 1] \cong [0,
1]$ as limit spaces or pseudotopological spaces. Consider the isomorphism
$\phi : [0, 1] \rightarrow [0, 1] \cup_{0 \sim 1} [0, 1]$ given by
\[ \phi (x, t) \mapsto \left\{\begin{array}{ll}
     q_0 (2 t ), & t \in \left[ 0, \frac{1}{2} \right)\\
     q_1 (2 t - 1), & t \in \left[ \frac{1}{2}, 1 \right],
   \end{array}\right. \]
where the maps $q_k : [0, 1] \rightarrow [0, 1] \cup_{1 \sim 0} [0, 1]$ are
the maps from the pushout square of $[0, 1] \cup_{1 \sim 0} [0, 1]$. To check
that $\phi$ is continuous, it suffices to check that $\phi$ is continuous at
the point $\frac{1}{2} \in [0, 1]$. Suppose that a filter $\lambda \rightarrow
\frac{1}{2}$ in the topological structure $([0, 1], \tau)$, and note that
$\phi (\lambda)$ is generated by the sets of the form $(t_0, 1] \cup_{1 \sim
0} [0, t_1) \subset [0, 1] \cup_{1 \sim 0} [0, 1]$. Sets of the form $(t_0,
1]$ generate the neighborhood filter $\mathcal{N}_1$ of the point $1$ in $[0,
1]$, and, similarly, sets of the form $[0, t_1)$ generate the neighborhood
filter $\mathcal{N}_0$ of $0$ in $[0, 1]$. The filter $\phi (\lambda)$
therefore contains the filter $q_0 (\mathcal{N}_1) \cap q_1 (\mathcal{N}_1)$,
which converges in quotient convergence structure of $[0, 1] \cup_{1 \sim 0}
[0, 1]$. It follows that $\phi (\lambda) \rightarrow \frac{1}{2}$, and
therefore $\phi$ is continuous.

Note that $\phi^{- 1}$ is continuous since it is the dashed arrow in the
pushout diagram
\begin{equation}
\begin{tikzcd}   
   I \ar[r,"i_0"] \ar[d,"i_1"] & I \ar[d,"q_1"] \ar[ddr,bend left=30,"\alpha_0"]& \\
   I \ar[r,"q_0"] \ar[drr,bend right=20,"\alpha_1"] & I \cup_{1 \sim 0} I \ar[dr,"\phi^{-1}",dashed] & \\
   & & I
   \end{tikzcd}
\end{equation}
where $\alpha_0, \alpha_1 : I \rightarrow I$ are given by $\alpha_k (t) =
\frac{t + k}{2}$ and $I = [0, 1]$.

Since $\phi, \phi^{- 1}$ are both bijective, we have that $I X \cong I X
\cup_X I X$ in limit spaces, as desired. The proof for pseudotopological
spaces is analogous.

\subsubsection{A Non-Example: General Convergence Spaces}In general
convergence spaces, it is no longer true that if the filters $\lambda_0
\rightarrow x$ and $\lambda_1 \rightarrow x,$ then $\lambda_0 \cap \lambda_1
\rightarrow x$, and therefore the spaces $[0, 1]$ and $[0, 1] \cup_{1 \sim 0}
[0, 1]$ are not isomorphic. See {\cite{Preuss_2002}}, Remark A.2.3 for a more
complete discussion of this point in the context of Kent convergence spaces.
General convergence spaces (and Kent convergence spaces) therefore do not
satisfy the homotopy gluing property.

\section{Convenient Cofibration Categories and Homotopy
Groups}\label{sec:I-category}

The general strategy for studying homotopy theory in a category is, first, to
specify a class of morphisms, often called weak equivalences, at which one
would like to localize the category, and, second, to endow the category with
enough structure in order to develop tools to probe the resulting homotopy
category. The structures which we will study in this article are model
categories {\cite{Quillen_1967,Hovey_1999,Hirschhorn_2019}} and the
cofibration category theory developed in {\cite{Baues_1989,Baues_1999}}. While
not as modern as more recent approaches to homotopy theory through infinity
categories {\cite{Riehl_Verity_2022,Cisinski_2019,Lurie_2009}}, they
nonetheless remain somewhat more accessible, and, once one has a model
category, an infinity category structure may also be constructed if desired.

In this section, we describe the axioms for cofibration categories and
$I$-categories from {\cite{Baues_1989}}, and we briefly indicate how these
structures may be used to construct homotopy groups. More details on these
constructions and the consequences of this theory may be found in
{\cite{Baues_1989,Baues_1999}}. After this brief introduction, we show that
convenient categories admit an $I$-category structure if they admit a cylinder
$(I, i_0, i_1, p)$ which satisfies the homotopy gluing property, and we deduce
from this the existence of a cofibration category structure on $\catname{Lim}$
and $\catname{PsTop}$. The resulting homotopy groups from the cofibration
category structure will then be used to construct a model category structure
on $\catname{PsTop}$ in Section \ref{sec:Model category}.

\subsection{Cofibration Categories and Homotopy}We begin with the definition
of of a cofibration category.

\begin{definition}[Cofibration Category]
  \label{def:Cofibration category}A \emph{cofibration category} is a tuple
  $(\mathcal{C}, \tmop{cof}, \tmop{we})$, where $\mathcal{C}$ is a category,
  and $\tmop{cof}$ and $\tmop{we}$, the \emph{cofibrations} and
  \emph{weak equivalences}, respectively, are special classes of morphisms
  of $\mathcal{C}$ which satisfy the following axioms
  \begin{enumerate}
    \item (Composition axiom) The isomorphisms in $\mathcal{C}$ are weak
    equivalences and also cofibrations. For two maps
    \[ A \xrightarrow{f} B \xrightarrow{g} C, \]
    if any two of $f, g$ and $g \circ f$ are weak equivalences, then so is the
    third. The composition of cofibrations is a cofibration.
    
    \item \label{item:Cof pushout axiom}(Pushout axiom) For a cofibration $i :
    B \rightarrowtail A$ and a map $f : B \rightarrow Y$, there exists a
    pushout in $\mathcal{C}$
    \begin{equation*}
    	\begin{tikzcd}
    		B \ar[d,"i",tail] \ar[r,"f"] & Y \ar[d,"\bar{i}",tail]\\
    		A \ar[r,"\bar{f}"] & A\cup_f Y
    		\end{tikzcd}
    \end{equation*}
    and $\bar{i}$ is a cofibration. Moreover,
    \begin{enumerate}
      \item if $f$ is a weak equivalence, then so is $\bar{f}$,
      \item if $i$ is a weak equivalence, then so is $\bar{i}$.
    \end{enumerate}
    \item \label{item:Factorization}(Factorization axiom) For a map $f : B
    \rightarrow Y$ in $\mathcal{C}$ there exists a commutative diagram
    \begin{equation*}
      \begin{tikzcd}  
       B \ar[rr,"f"] \ar[rd,tail,"i"] & & Y\\
       & A \ar[ru,"g","\sim"']
      \end{tikzcd}
  	\end{equation*}
    where $i$ is a cofibration and $g$ is a weak equivalence.
    
    \item (Axiom on fibrant models) For each object $X$ in $\mathcal{C}$ there
    is a trivial cofibration $X \rightarrowtail R X$ (i.e. a cofibration which
    is also a weak equivalence) where $R X$ is \emph{fibrant} in
    $\mathcal{C}$, i.e. each trivial cofibration $i : R X \rightarrowtail Q$
    admits a retraction $r : Q \rightarrow R X, r i = 1_{R X}$. 
  \end{enumerate}
\end{definition}

We now summarize the construction of homotopy groups in a general cofibration
category. We will discuss the specific cases of $\catname{Lim}$ and
$\catname{PsTop}$ in Section \ref{subsubsec:PsLim homotopy}. The following
discussion is summarized from \cite{Baues_1989}, II, Sections 1-6. We assume
throughout that $\mathcal{C}$ is a cofibration category with an initial object
$\ast$ such that $\tmop{Hom} (A, \ast)$ is non-empty for any
$\mathcal{C}$-object $A$. Note that the latter condition is automatically
satisfied if $\ast$ is also a terminal object.
\begin{definition}
  Relative Cylinders. Let $i : B \rightarrowtail A$ be a cofibration. Then, by
  the Pushout Axiom (Definition \ref{def:Cofibration category}, Item
  \ref{item:Cof pushout axiom}) there exists the pushout
  \begin{equation*}
     \begin{tikzcd}
     B \ar[d,tail]\ar[r,tail] & A \ar[d,tail]\ar[ddr,"1_A",bend left=30]\\   
     A \ar[r,tail] \ar[drr,"1_A",bend right=30]& A \cup\_B A \ar[dr,dashed,"\phi"]\\
     & & A
     \end{tikzcd}
 \end{equation*}
  We call $\phi = (1_A, 1_A)$ the \emph{folding map}. By the Factorization
  Axiom (Definition \ref{def:Cofibration category}, Item
  \ref{item:Factorization}), there is a factorization of $\phi$
  \[ A \cup_B A {{\xrightarrow{i} } }  Z \xrightarrow[\sim]{p} A \]
  where $i$ is a cofibration and $p$ is a weak equivalence. We call the triple
  $(Z, i, p)$ a \emph{relative cylinder on $B \rightarrowtail A$}. We
  often write $I_B A$ or $Z$ for the relative cylinder $(Z, i, p)$.
\end{definition}

Note the difference between relative cylinders and the cylinders of Definition
\ref{def:Cylinder}.

\begin{definition}
  Homotopy. Let $X$ be a fibrant object in a cofibration category
  $\mathcal{C}$, and let $Y \subset X$ be a cofibration. We say that two maps
  $f, g : A \rightarrow X$ are \emph{homotopy relative $B$}, denoted $f
  \simeq g \tmop{rel} B$, iff there is a commutative diagram
  \begin{equation*}     
     \begin{tikzcd}[column sep=1.5em]   
     A \cup_B A \ar[dr,"(f{,}g)"']\ar[rr,tail,"i"] & & Z \ar[dl,"H"] \\
     & X
     \end{tikzcd}
 \end{equation*}
  where $Z$ is a relative cylinder on $B \rightarrowtail A$. We call $H : Z
  \rightarrow X$ a \emph{homotopy rel $B$ from f to g}.
\end{definition}

\begin{definition}
  Let $Y \subset X$ be a cofibration in a cofibration category $\mathcal{C}$.
  Let $u : Y \rightarrow U$ be a map and let $\tmop{Hom} (X, U)^u$ denote the
  set of all \emph{extensions of $u$}, i.e. the maps \ $f : X \rightarrow
  U$ in $\mathcal{C}$ such that $f \mid_Y = u.$
\end{definition}

\begin{proposition}[{\cite{Baues_1989}}, Proposition II.2.2]
  Suppose that $U$ is a fibrant object in a cofibration category
  $\mathcal{C}$, and let $Y \subset X$ be a cofibration. Then all of the
  cylinders on $Y \subset X$ define the same homotopy relation relative $Y$,
  denoted $\simeq$ $\tmop{rel}$ $Y$, on the set $\tmop{Hom} (X, U)^u$, and,
  furthermore, the homotopy relation relative $Y$ is an equivalence relation.
\end{proposition}

\begin{definition}
  Torus. We define the \emph{torus $\Sigma_Y X$} for $Y \subset X$ by the
  pushout diagram
  \begin{equation*}
     \begin{tikzcd}
     X \cup\_Y X \ar[d,tail]
     \ar[r,"(1{,}1)"] & X \ar[d,tail,"i"]
     \ar[rdd,bend left=40,"1"]\\
     I_Y X \ar[r,"\tau"] \ar[drr,bend right=30,"p"] & \Sigma_Y X \ar[dr,dashed,"r"]\\
     & & X
     \end{tikzcd}
 \end{equation*}
  where $p$ is the map in the definition of the cylinder $(I_Y X, i, p)$ and
  the existence of $r$ is guaranteed by the universal property of the pushout.
  
\end{definition}

\begin{definition}
  Based object. A \emph{based object} in a cofibration category
  $\mathcal{C}$ is a pair $(X, o_X)$, where $X$ is a cofibrant object (i.e.
  $\ast \rightarrow X$ is a cofibration) and where $o \assign o_X : X
  \rightarrow \ast$ is a map from $X$ to the initial object. We call $o = o_X$
  the \emph{trivial map on $X$}. We say that a map $f : (X, o_X)
  \rightarrow (Y, o_Y)$ between based objects is \emph{based} iff $o_Y
  \circ f = o_X$. (Note that in general there may be many maps $X \rightarrow
  \ast$ to the initial object.)
\end{definition}

\begin{definition}
  Suspension of a based object. For a based object $A$ in a cofibration
  category $\mathcal{C}$ the suspension $\Sigma A$ is the based object which
  is defined by the pushout diagram
	\begin{equation*}
     \begin{tikzcd}     
     A \vee A \ar[d,"(i_0{,}i_1)"',tail]\ar[r,"(1{,}1)"] \ar[dr,phantom,"push"]& A \ar[d,tail]\ar[dr,"push",phantom]\ar[r,"o_A"]
     & \ast \ar[d,tail] \ar[dd,bend left=40]\\
     I_\ast A \ar[dr,"p"]\ar[r] &
     \Sigma_\ast A \ar[d,"r"]\ar[r,"\sigma"] & \Sigma A \ar[d,"o_{\Sigma A}"']\\
     & A \ar[r,"o_A"] & \ast
     \end{tikzcd}
  \end{equation*}
  where $\Sigma_{\ast} A$ is the \emph{torus on $\ast \rightarrowtail A$}.
  Note that the existence of the map $o_{\Sigma A} : \Sigma A \rightarrow
  \ast$ exists by the universal property of the pushout. Also note that, a
  priori, $\Sigma A$ depends on the choice of trivial map $o_A : A \rightarrow
  \ast$. However, in the cases of interest in this article, $\ast$ will be a
  terminal object as well as an initial object, so $\Sigma A$ will be
  canonical.
\end{definition}

\begin{definition}
  For a based object $A$ in a cofibration category $\mathcal{C}$, we define
  $\pi_n^A (U) \assign [\Sigma^n A, U]$. These are groups for $n \geq 1$ by
  {\cite{Baues_1989}}, II.5.13 (As noted in the discussion in
  {\cite{Baues_1989}}, II.6.9-11.)
\end{definition}

\begin{example}
  Homotopy groups of topological spaces. Let $\mathcal{C} =
  \catname{Top^{\ast}}$, the category of pointed topological spaces
  (with initial object $\ast$), and let $A = S^0$, the $0$-dimensional sphere.
  Then $\Sigma^n S^0$ is the $n$-th based suspension of $S^0$, which is
  homeomorphic to $S^n$, and $\pi_n (X) = \pi_n^{S^0} (X) = [\Sigma^n S^0,
  X]$, which are the usual homotopy groups in $\catname{Top^{\ast}}$.
\end{example}

\subsection{$I$-categories}\label{subsec:I-categories}

We will not work directly with the axioms of a cofibration category, but
instead with an auxiliary structure called an \emph{$I$-category}, which
is an axiomatization of the case when the category in question has a canonical
cylinder, as is the case for $\catname{Lim}$ and $\catname{PsTop}$, as well as
for convenient categories in general. An $I$-category is defined as follows.
Note that we have slightly reorganized the axioms in {\cite{Baues_1989}},
permuting the last two axioms and incorporating much of the cylinder axiom in
the definition of the cylinder $(I, i_0, i_1, p)$.

\begin{definition}
  \label{def:I-category}An \emph{$I$-category} is a tuple $(\mathcal{C},
  \tmop{cof}, (I, i_0, i_1, p), {\varnothing})$ where $\mathcal{C}$ is a category,
  $\tmop{cof}$ is a collection of $\mathcal{C}$-morphisms called cofibrations,
  $(I, i_0, i_1, p)$ is a cylinder (Definition \ref{def:Cylinder}) on
  $\mathcal{C}$, and ${\varnothing}$ is the initial object in $\mathcal{C}$. The
  structure must satisfy the following axioms:
  \begin{enumerate}
    \item \label{item:Cylinder axiom}(Cylinder axiom) $I{\varnothing}={\varnothing}$.
    
    \item \label{item:Pushout axiom}(Pushout axiom) For a cofibration $i : B
    \rightarrowtail A$ and a morphism $f : B \rightarrow X$ there exists the
    pushout
    \[
       \begin{tikzcd}
       B \ar[r,"f"] \ar[d,rightarrowtail,"i"]& X
       \ar[d,rightarrowtail,"\bar{i}"]\\
       A \ar[r] & A\cup_B X
       \end{tikzcd} \]
    and $\bar{i}$ is a cofibration. Moreover, the functor $I$ carries pushouts
    to pushouts, i.e. $I (A \cup_B X) = I A \cup_{I B} I X.$
    
    \item \label{item:Cofibration axiom}(Cofibration axiom) The class
    $\tmop{cof}$ of cofibrations satisfy the following:
    \begin{enumerate}
      \item \label{item:Cofibration a}Each isomorphism is a cofibration
      
      \item For every $\mathcal{C}$-object $X$, the morphism ${\varnothing} \rightarrow
      X$ is a cofibration
      
      \item The composition of cofibrations is a cofibrations
      
      \item \label{item:HEP}A cofibration $i : B \rightarrowtail A$ satisfies
      the following homotopy extension property. For each commutative solid
      arrow diagram
      
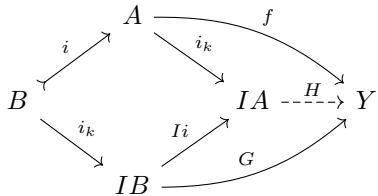
\[
         \begin{tikzcd}
         & A \ar[dr,"i_k"]\ar[drr,bend
         left=20,"f"] & \\
         B \ar[ur,rightarrowtail,"i"]\ar[dr,"i_k"] & & IA \ar[r,dashed,"H"]&
         Y\\
         & IB \ar[ur,"Ii"]\ar[urr,bend right=20,"G"]
         \end{tikzcd} \]
      in $C$, there exists an $H : I A \rightarrow Y$ which makes the entire
      diagram commute.
    \end{enumerate}
    \item (Interchange axiom) For all $\mathcal{C}$-objects $X$, there exists
    a morphism $T : I I X \rightarrow I I X$ with $T i_k = I i_k$ and $T I
    (i_k) = i_k$ for $k = 0, 1$. $T$ is called the \emph{interchange map.}
    
    \item \label{item:Relative cylinder axiom}(Relative Cylinder axiom) For a
    cofibration $i : B \rightarrowtail A$, the map $j$ defined by the pushout
    diagram
    
    is a cofibration.
  \end{enumerate}
  We will often abbreviate $(\mathcal{C}, (I, i_0, i_1, p), \tmop{cof}, {\varnothing})$
  by $(\mathcal{C}, I, \tmop{cof}, {\varnothing})$ when the remainder of the cylinder
  is clear from context.
\end{definition}

The main interest of an $I$-category structure is that, when a category
$\mathcal{C}$ does have a canonical cylinder functor, it is often easier to
verify the $I$-category axioms than it is to directly verify the cofibration
category axioms. However, the following theorem from {\cite{Baues_1999}} shows
that an $I$-category structure implies the existence of a cofibration
category. We begin by defining homotopy equivalence in an $I$-category.

\begin{definition}
  \label{def:Homotopy I-cat}(Homotopy and homotopy equivalence in an
  $I$-category)Let $(\mathcal{C}, I, \tmop{cof}, {\varnothing})$ be an $I$-category. We
  say that two maps $f, g : X \rightarrow Y$ are \emph{homotopic}, written
  $f \simeq g$, if there exists a morphism $H : \tmop{IX} \rightarrow Y$ such
  that $H i_0 = f$ and $H i_1 = g$. We call $H$ the \emph{homotopy from
  $f$ to $g$}. We say that a map $f : X \rightarrow Y$ is a \emph{homotopy
  equivalence} if there is a map $g : Y \rightarrow X$ such that $f g \simeq
  1_Y$ and $g f \simeq 1_X$.
\end{definition}

\begin{theorem}[{\cite{Baues_1999}}, Theorem 3.3]
  \label{thm:I-cat is cof cat}Let $(\mathcal{C}, \tmop{cof}, I, {\varnothing})$ be an
  $I$-category. Then $\mathcal{C}$ is a cofibration category with the
  following structure. Cofibrations are those of the $I$-category structure,
  and weak equivalences are the homotopy equivalences from Definition
  \ref{def:Homotopy I-cat}.
\end{theorem}

Note that in this structure, all objects are fibrant and cofibrant in
$\mathcal{C}$.
\subsection{An $I$-category structure on convenient categories with the
homotopy gluing property}

We now show that any convenient category $\mathcal{C}$ admits a cylinder $(I,
i_0, i_1, p)$, and then that if the cylinder satisfies the homotopy gluing
property, then $\mathcal{C}$ admits a cofibration category structure.

\begin{proposition}
  Let $\mathcal{C}$ be a category which contains the topological interval $[0,
  1]$ and the one point set $\{ \ast \}$ as objects, such that the set maps
  $\ast \mapsto 0 \in [0, 1]$ and $\ast \mapsto 1 \in [0, 1]$ are elements of
  $\tmop{Hom}_{\mathcal{C}} (\{ \ast \}, [0, 1])$. Then there exists a
  cylinder $(I, i_0, i_1, p)$ on $\mathcal{C}$.
\end{proposition}

\begin{proof}
  We define the functor $I - : \mathcal{C} \rightarrow \mathcal{C}$ to be the
  product functor $- \times I$, and for each $\mathcal{C}$-object $X$, we
  define the $i_{k, X}$ to be the unique dashed arrow in the diagram
  \begin{equation}\label{eq:Pullback diagram}
    \begin{tikzcd}
    X \ar[d]\ar[drr,"Id_X",bend
    left=20]\ar[dr,dashed,"i_{k,X}" description] & & \\
    \ast \ar[dr,"\iota_k"']& IX
    \ar[d]\ar[r] & X \ar[d]\\
    & {[}0,1{]}\ar[r] & \ast
    \end{tikzcd}
  \end{equation}
  The unicity of the dashed arrow in the diagram
  \begin{equation*}	 
     \begin{tikzcd}
     X \ar[d]\ar[dr,dashed] \ar[drr,"f"] & &\\
     \ast \ar[dr,"\iota_k"] & IY \ar[r]\ar[d] & Y \ar[d] \\
     & {[}0,1{]} \ar[r] & \ast
    \end{tikzcd}
  \end{equation*}
  implies the commutativity of the diagram
  \begin{equation*}
     \begin{tikzcd}   
     X \ar[d,"i_{k,X}"]\ar[r,"f"] & Y\ar[d,"i_{k,Y}"]\\
     IX \ar[r,"If"]& IY
     \end{tikzcd}
   \end{equation*}
  and therefore $i_0, i_1 : \tmop{Id}_{\mathcal{C}} \rightarrow I$ are natural
  transformations.
  Similarly, for each
  $\mathcal{C}$-object $X$ we define $p_X : \tmop{IX} \rightarrow X$ to be the
  map $I X \rightarrow X$ from the pullback diagram (\ref{eq:Pullback
  diagram}) of $I X.$ We see from this diagram that the composition
  \[ p \circ i_k (X) : X \rightarrow I X \rightarrow X \]
  is the identity of $X$ by construction, for all $X \in \mathcal{C}$ and $k =
  0, 1$. Therefore, $(I, i_0, i_1, p)$ is a cylinder in $\mathcal{C}$.
\end{proof}

The following corollary is immediate.

\begin{corollary}
  \label{prop:Convenient cylinder}A convenient category $\mathcal{C}$ has a
  cylinder $(I, i_0, i_1, p)$ defined by
  \begin{eqnarray*}
    I X & \assign & X \times [0, 1]\\
    I f (x, t) & \assign & f \times 1_{[0, 1]}\\
    i_0 X & : = & 1_X \times \{ 0 \}\\
    i_1 X & \assign & 1_X \times \{ 1 \}\\
    p (I X) & \assign & p_X (I X)
  \end{eqnarray*}
  where $p_X : X \times [0, 1] \rightarrow X$ is the projection onto $X$ in
  the universal property of the product.
\end{corollary}

We now define the cofibrations in a convenient category.

\begin{definition}
  \label{def:Cofibrations}Let $\mathcal{C}$ be an convenient category. Then we
  define $\tmop{cof}$ to be the collection of $\mathcal{C}$-morphisms which
  satisfy the homotopy extension property in the cofibration axiom, i.e. Item
  \ref{item:Cofibration axiom}\ref{item:HEP} of Definition
  \ref{def:I-category}.
\end{definition}

We now check the axioms of an $I$-category for a convenient category  C  whose
cylinder $(I, i_0, i_1, p)$ satisfies the homotopy gluing property, with the
cofibrations $\tmop{cof}$ defined as in Definition \ref{def:Cofibrations}
above. We take the initial object ${\varnothing}$ to be the empty set.

\begin{proposition}
  Let $\mathcal{C}$ be a convenient category, the cylinder $(I, i_0, i_1, p)$
  be defined as in Proposition \ref{prop:Convenient cylinder}, and let
  $\tmop{cof}$ be the class of cofibrations defined in Definition
  \ref{def:Cofibrations}. Then $(\mathcal{C}, I, \tmop{cof}, {\varnothing})$ satisfies
  the Pushout Axiom.
\end{proposition}

\begin{proof}
  Since $\mathcal{C}$ is Cartesian closed by hypothesis, pushouts exist and
  $I$ preserves pushouts. It remains to verify that $\bar{i}$ is a
  cofibration. That is, given the commutative diagram of solid arrows
  \begin{equation*}	
     \begin{tikzcd}
     & A \cup_B Y \ar[dr, "i_0"] \ar[drr, bend left=20, "f"] \\
     Y \ar[dr,"i_0"] \ar[ur,"\bar{i}"]& & I(A \cup_B Y)\ar[r, dashed,
     "H"] & Z\\ 
     & IY \ar[ur, "I\bar{i}"] \ar[urr, bend right=20, "G"]
     \end{tikzcd}
  \end{equation*}
  we wish to show that there is a morphism $H : I (A \cup_B Y) \rightarrow Z$
  which makes the entire diagram commute.
  
  First, we note that $I (A \cup_B Y) = I A \cup_{I B} I Y$. Second, $A
  \cup_B Y$ is a pushout, and therefore for any commutative diagram of solid
  arrows
	\begin{equation*}
     \begin{tikzcd}     
     & A \ar[dr,"q\_A"] \ar[drr,bend left=20,"f\circ q A"] & & \\
     B \ar[dr,"\beta_B"] \ar[ur,"i"] & & A\cup\_B Y \ar[r,dashed,"f"] & Z\\
     & Y \ar[ur,"\bar{i}"] \ar[urr,bend right=20,"f\circ \bar{i}"]
     \end{tikzcd}
    \end{equation*}
  Since $I B \cup_{I B} I Y$ is also a pushout, we also have
  \[ 
     \begin{tikzcd}  
     & IA \ar[dr,"Iq\_A"]\ar[drr,bend left=20,"H"] & &\\ 
     IB \ar[dr,"I\beta_Y"]\ar[ur,"Ii"] & & IA\cup_{IB} IY \ar[r,dashed,"F"]& Z\\
     & IY \ar[ur,"I\bar{i}"] \ar[urr,bend right=20,"G"]
     \end{tikzcd} \]
  Since $i : B \rightarrow A$ is a cofibration, we know that, for any
  $\mathcal{C}$-object $Z$
  \[
     \begin{tikzcd}     
     &A \ar[dr, "i_{k,A}"] \ar[drr, bend left=20, "f\circ q_A"]\\
     B\ar[dr,"i_{k,B}"] \ar[ur,"i"]& &
     IA\ar[r, dashed, "H"] & Z\\
     & IB \ar[ur, "Ii"] \ar[urr, bend right=19,"G\circ I\beta_Y"]
     \end{tikzcd} \]
  Putting these all together, we have the diagram
  \begin{equation}\label{diag:Cofibration}
    \begin{tikzcd}
    & A \ar[dr,"q_A"]\ar[ddrrr,"f\circ q_A",bend left] & \\
    B \ar[dr,"\beta\_B"]\ar[ur,"i"]& & A \cup_B Y \ar[dr,"i_{k,A
    	\cup_B Y}"] \ar[drr, bend left=20, "f"]\\
    & Y\ar[dr,"i_{k,Y}"] \ar[ur,"\bar{i}"]& & I(A \cup_B Y)\ar[r, dashed, "F"] &[2.5em] Z\\
    & & IY \ar[ur, "I\bar{i}"] \ar[urr, bend right=20, "G"]
    \end{tikzcd}
  \end{equation}
  where the solid-arrow diagram commutes. Finally, since the $i_k$ are natural
  transformations, we have
  \[ \begin{tikzcd}
     A \ar[d,"i_{k,A}"]\ar[r,"q_A"] &
     A\cup_B Y \ar[d,"i_{k,A \cup\_B Y}"]\\
     IA \ar[r,"Iq_A"] & I(A\cup_B)
     \end{tikzcd} \]
  From these last two diagrams, we conclude that
  \begin{eqnarray*}
    f \circ q_A & = & H \circ i_{k, A}\\
    & = & F \circ I q_A \circ i_{k, A}\\
    & = & F \circ i_{k, A \cup_B Y} \circ q_A .
  \end{eqnarray*}
  We also have
  \begin{eqnarray*}
    f \circ \bar{i} & = & G \circ i_{k, Y}\\
    & = & F \circ I \bar{i} \circ i_{k, Y}\\
    & = & F \circ i_{k, A \cup_B Y} \circ \bar{i} .
  \end{eqnarray*}
  Since both $f$ and $F \circ i_{k, A \cup_B Y}$ satisfy the pushout property
  of Diagram \ref{diag:Cofibration}, we conclude that $f = F \circ i_{k, A
  \cup_B Y}$ by the unicity of $f$. It now follows that $F : I (A \cup_B Y)
  \rightarrow Z$ makes the complete Diagram \ref{diag:Cofibration} commute,
  and since $Z$ was an arbitrary $\mathcal{C}$-object, we conclude that
  $\bar{i}$ is a cofibration.
\end{proof}

\begin{proposition}
  Let $\mathcal{C}$ be a convenient category, the cylinder $(I, i_0, i_1, p)$
  be defined as in Proposition \ref{prop:Convenient cylinder}, and let
  $\tmop{cof}$ be the class of cofibrations defined in Definition
  \ref{def:Cofibrations}. Then $(\mathcal{C}, I, \tmop{cof}, {\varnothing})$ satisfies
  the Cofibration Axiom.
\end{proposition}

\begin{proof}
  We prove Items \ref{item:Cofibration a}-\ref{item:HEP} in the cofibration
  axiom one by one.
  \begin{enumeratealpha}
    \item Each isomorphism is a cofibration: Let $i : B \rightarrow A$ be an
    isomorphism. Then $I i : I B \rightarrow I A$ is an isomorphism by
    construction, and for any diagram of solid arrows
    
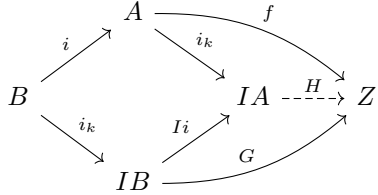
\[ 
       \begin{tikzcd}& A \ar[dr,"i_k"]\ar[drr,bend left=20,"f"] \\
       B \ar[ur,"i"]\ar[dr,"i_k"] & & IA \ar[r,dashed,"H"]& Z\\
       & IB \ar[ur,"Ii"]\ar[urr,bend right=20,"G"]
       \end{tikzcd} \]
    we define $H : I A \rightarrow Z$ by $H \assign G \circ (I i)^{- 1}$, and
    $H$ makes the full diagram commute. Since $Z$ is arbitrary, $i$ is a
    cofibration.
    
    \item For every $\mathcal{C}$-object $X$, the morphism ${\varnothing} \rightarrow
    X$ is a cofibration: Since $I{\varnothing}={\varnothing}$ by hypothesis, for each diagram of
    solid arrows
    
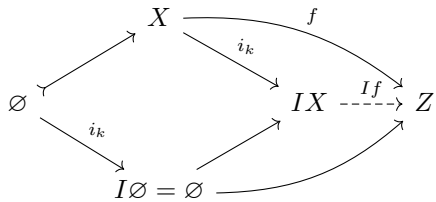
\[ 
       \begin{tikzcd}& X \ar[dr,"i_k"]\ar[drr,bend left=20,"f"] &\\
       \varnothing \ar[ur,rightarrowtail]\ar[dr,"i_k"] & & IX \ar[r,dashed,"If"]&
       Z\\
       & I\varnothing=\varnothing \ar[ur]\ar[urr,bend right=20]
       \end{tikzcd} \]
    the map $I f : I X \rightarrow Z$ makes the full diagram commute.
    
    \item The composition of cofibrations is a cofibration: Let $i : B
    \rightarrowtail A$ and $j : C \rightarrowtail B$ be cofibrations, and
    consider the map $i \circ j$. We wish to show that, for any diagram of
    solid arrows
    
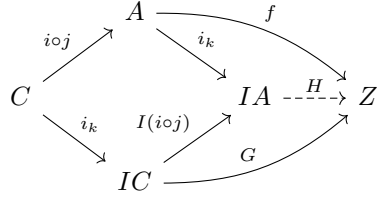
\[
       \begin{tikzcd}
       & A \ar[rrd,bend left=20,"f"] \ar[rd,"i_k"]& \\ 
       C \ar[ru,"i\circ j"] \ar[rd,"i_k"] & & IA \ar[r,"H",dashed] & Z\\
       & IC \ar[ru,"I(i\circ j)"] \ar[rru,bend right=20,"G"]&
       \end{tikzcd}
   \]
    there exists a dashed arrow $H : I A \rightarrow Z$ making the full
    diagram commute. We first consider the following diagram
    
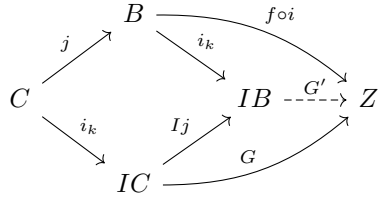
\[ 
       \begin{tikzcd}
       & B \ar[rrd,bend left=20,"f\circ i"]
       \ar[rd,"i_k"]& \\
       C \ar[ru,"j"] \ar[rd,"i_k"] & & IB
       \ar[r,"G'",dashed] & Z\\
       & IC \ar[ru,"Ij"] \ar[rru,bend
       right=20,"G"]&
       \end{tikzcd} \]
    Since $j : C \rightarrowtail B$ is a cofibration, there exists $H' : I B
    \rightarrow Z$ making the above diagram commute. However, since $i : B
    \rightarrowtail A$ is a cofibration, there exists $H : I A \rightarrow Z$
    making the diagram
    
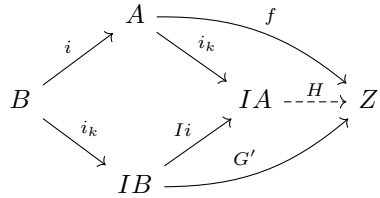
\[
       \begin{tikzcd}       
       & A \ar[rrd,bend left=20,"f"] \ar[rd,"i_k"]& \\
       B \ar[ru,"i"] \ar[rd,"i_k"] & & IA \ar[r,"H",dashed] & Z\\
       & IB \ar[ru,"Ii"] \ar[rru,bend right=20,"G'"] &
       \end{tikzcd}
    \]
    commute. Putting these together, we have
    
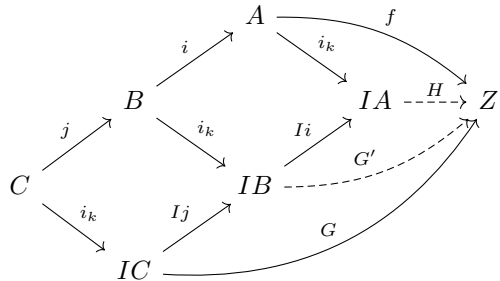
\[
       \begin{tikzcd}
       & & A \ar[rrd,bend left=20,"f"] \ar[rd,"i_k"]& \\
       & B \ar[ru,"i"] \ar[rd,"i_k"] & & IA
       \ar[r,"H",dashed]& Z\\
       C \ar[ru,"j"] \ar[rd,"i_k"] & & IB
       \ar[ru,"Ii"] \ar[rru,"G'",dashed,bend right=20]\\
       & IC \ar[ru,"Ij"] \ar[rrruu,bend
       right=30,"G"]&
       \end{tikzcd} \]
    as desired.
    
    \item A cofibration $i : B \rightarrowtail A$ in $(\mathcal{C}, I,
    \tmop{cof}, {\varnothing})$ satisfies the homotopy extension property by
    definition.
  \end{enumeratealpha}
  
\end{proof}

For the proof of the Relative Cylinder Axiom, we will require the following a
preparatory lemma.

\begin{lemma}
  A map $f : B \rightarrow A$ is a cofibration as defined in Definition
  \ref{def:Cofibrations} iff there exists a retract $r_k : I A \rightarrow I B
  \cup_{i_k} A$ of the map $j_k : I B \cup_{i_k} A \rightarrow I A$ in the
  pushout diagram
  
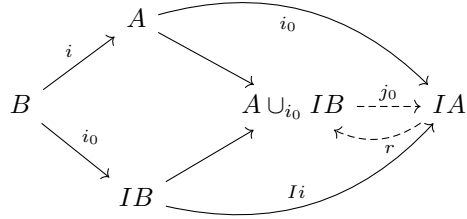
\[
     \begin{tikzcd}
     & A \ar[dr] \ar[drr,"i_0"',bend
     left=30]\\
     B \ar[dr, "i_0"] \ar[ur,"i"] & & A
     \cup_{i_0} IB \ar[r,"j_0",dashed] & IA
     \ar[l, bend left, dashed, "r"] \\  
     & IB \ar[ur] \ar[urr,"Ii",bend right=30]  
     \end{tikzcd} \]
\end{lemma}

\begin{proof}
  First, suppose that $i : B \rightarrow A$ is a cofibration. Then for any
  $Y$ and maps $f : X \rightarrow Y$ and $G : I A \rightarrow Y$, there exists
  a map $H : I X \rightarrow Y$ that makes the diagram
  \[ 
     \begin{tikzcd} 
     & A \ar[dr, "i_0"] \ar[drr, bend left=20,
     "f"]\\ 
     B \ar[dr,"i_0"] \ar[ur,"i"]& &
     IA\ar[r, dashed, "H"] & Y\\
     & IB \ar[ur, "Ii"] \ar[urr, bend right=20,"G"]
     \end{tikzcd} \]
  commute. Let $Y \assign A \cup_{i_k} I B$, and let $f : A \rightarrow A
  \cup_{i_k} I B$ and $G : I B \rightarrow A \cup_{i_k} I B$ be the standard
  inclusions for the pushout. Then the map $r = H : I A \rightarrow A
  \cup_{i_k} I B$ is the desired retract, by the commutativity of the diagram.
  
  Now suppose such a retract $r : I A \rightarrow A \cup_{i_k} I B$ of $j_k$
  exists. Then, for any $\mathcal{C}$-object $Y$ and maps $f : A \rightarrow
  Y$ and $G : I B \rightarrow Y$ such that the solid arrow diagram
  
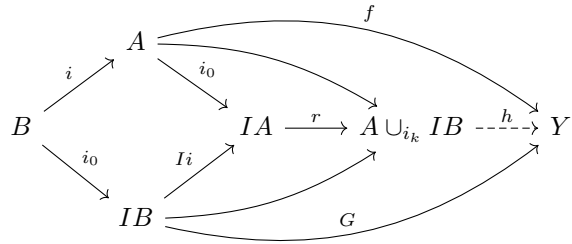
\[ 
     \begin{tikzcd}
     & A \ar[dr, "i_0"] \ar[drr,bend left=15]\ar[drrr, bend left=25,
     "f"]\\
     B \ar[dr,"i_0"] \ar[ur,"i"]& & IA
     \ar[r,"r"]& A\cup_{i_k} IB
     \ar[r, dashed, "h"] & Y\\
     & IB \ar[ur, "Ii"] \ar[urr,bend right=15]
     \ar[urrr, bend right=25, "G"]
     \end{tikzcd} \]
  exists, then there exists a map $h : A \cup_{i_k} I B \rightarrow Y$ making
  the entire diagram commute, Defining $H \assign h \circ r$, we see that $i$
  is a cofibration, as desired.
\end{proof}

\begin{proposition}
  Let $\mathcal{C}$ be a convenient category, the cylinder $(I, i_0, i_1, p)$
  be defined as in Proposition \ref{prop:Convenient cylinder}, and let
  $\tmop{cof}$ be the class of cofibrations defined in Definition
  \ref{def:Cofibrations}. Suppose, in addition, that the cylinder $(I, i_0,
  i_1, p)$ satisfies the homotopy gluing property in $\mathcal{C}$. Then
  $(\mathcal{C}, I, \tmop{cof}, {\varnothing})$ satisfies the Relative Cylinder Axiom.
\end{proposition}

\begin{proof}
  We follow the proof in {\cite{Baues_1999}} for topological spaces. For the
  map $j : A \cup I B \cup A \rightarrow I A$ to be a cofibration, for any
  $\mathcal{C}$-object $X$ and any diagram of solid arrows
  \[
     \begin{tikzcd}
     & IA \ar[dr,"i_k"] \ar[drr,"f",bend left=20]& & \\
     A\cup IB \cup A \ar[ur,"j"] \ar[dr,"i_k"]& & IIA \ar[r,"H",dashed]& X\\
     & I(A \cup IB \cup A)\ar[ur,"Ij"] \ar[urr,"G",bend right=20]& &
     \end{tikzcd}\]
  there must exist a dashed arrow $H : I I A \rightarrow X$ making the full
  diagram commute. \ Let $\gamma : I \times I \rightarrow \mathbb{R}$ be the
  map $\gamma (x, t) \assign \max (2 \| x \|, 2 - t)$, and define the
  homeomorphism $\alpha : I \times I \rightarrow I \times I$ by $\alpha (x, t)
  : = (\gamma (x, y)^{- 1} (1 + t) x, 2 - \gamma (x, t))$. Consider the
  diagram
  \[
     \begin{tikzcd}     
     (IA) \cup_{i_0} I(A \cup_{i_0 B} IB
     \cup_{i_1 B} A) \ar[d,"\beta"] \ar[r, "j_0"] & IIA \ar[d,"\alpha \times 1_A"]\\  
     I(A \cup_B IB) \ar[r,"TI\iota"] & IIA
    \end{tikzcd}
\]
  where $\beta$ is the restriction of $\alpha \times 1_A$ which makes the
  diagram commute, $\iota : A \cup_B I B \rightarrow I A$ is the inclusion,
  and $T : I I A \rightarrow I I A$ is the map from the Interchange Axiom.
  Since $\mathcal{C}$ satisfies the Homotopy Gluing Property, $\beta$ is a
  homeomorphism onto $I (A \cup_B I B)$, and therefore the inverse map
  $\beta^{- 1}$ exists. (Note that without the homotopy gluing property, the
  inverse $\beta^{- 1}$ may not be continuous.) Defining $\bar{r} \assign
  \beta^{- 1} (I r) T (\alpha \times 1_A)$ where $r : I A \rightarrow A \cup_B
  I B$ is the retract of $\iota$, which exists because $i : B \rightarrow A$
  is a cofibration. We see that $\bar{r}$ this is the desired retract of
  $j_0$, and therefore $j$ is a cofibration, as
  desired.
\end{proof}

\begin{proposition}
  Let $\mathcal{C}$ be a convenient category, the cylinder $(I, i_0, i_1, p)$
  be defined as in Proposition \ref{prop:Convenient cylinder}, and let
  $\tmop{cof}$ be the class of cofibrations defined in Definition
  \ref{def:Cofibrations}. Then $(\mathcal{C}, I, \tmop{cof}, {\varnothing})$ satisfies
  the Interchange Axiom.
\end{proposition}

\begin{proof}
  Let $I_0$ and $I_1$ be distinct copies of $I$. Then the limits $I_0 \times
  I_1$ and $I_1 \times I_0$ exist and there exists a $\mathcal{C}$-morphism
  $T_I : I_1 \times I_0 \rightarrow I_0 \times I_1$ from the properties of the
  limits $I_0 \times I_1$ and $I_1 \times I_0$. We let $T : I I A \rightarrow
  I I A $ be the morphism $T : = T_I \times \tmop{Id}_A$. Then
  \[ T (i_k (I A)) = T (\{ k \} \times I \times A) = T_I \times \tmop{Id}_A
     (\{ k \} \times I \times A) = (I \times \{ k \} \times A) = I (i_k A) \]
  for $k \in \{ 0, 1 \}$. We also have
  \[ T (I (i_k (A)) = T (I \times \{ k \} \times A) = T_I \times \tmop{Id}_A
     (I \times \{ k \} \times A) = \{ k \} \times I \times A = i_k (I A)
     \nobracket \]
  as desired.
\end{proof}

Putting these together, we have shown

\begin{theorem}
  \label{thm:I-category}An admissible category $(\mathcal{C}, I, \tmop{cof},
  {\varnothing})$ is an $I$-category, with $I, \tmop{cof}$, and ${\varnothing}$ defined as
  above.
\end{theorem}

The following corollaries are immediate.

\begin{corollary}
  \label{cor:PsTop I-category}$\left( \catname{PsTop}, I, \tmop{cof},
  {\varnothing} \right)$ is an $I$-category.
\end{corollary}

\begin{corollary}
  $\label{cor:Lim I-category} \left( \catname{Lim}, I, \tmop{cof}, {\varnothing}
  \right)$ is an $I$-category.
\end{corollary}

Combining Corollaries \ref{cor:PsTop I-category} and \ref{cor:Lim I-category}
with Theorem \ref{thm:I-cat is cof cat}, we have

\begin{corollary}
  $\left( \catname{PsTop}, \tmop{cof}, \tmop{we} \right)$ and $\left(
  \catname{Lim}, \tmop{cof}, \tmop{we} \right)$ are cofibration
  categories, where $\tmop{cof}$ is defined in Definition
  \ref{def:Cofibrations} and the weak equivalences are as in Theorem
  \ref{thm:I-cat is cof cat}.
\end{corollary}

\subsubsection{On the homotopy groups and suspensions of spheres in
$\catname{PsTop}$ and $\catname{Lim}$}\label{subsubsec:PsLim
homotopy}

We now show that a topological $n$-sphere can be obtained by the
pseudotopological based suspension of a topological $(n - 1)$-sphere. It
follows that the homotopy groups constructed from the cofibration category
structure on pseudotopological spaces are, in fact, the homotopy classes of
maps from the $n$-spheres, just as in the topological case. We begin with the
following lemma about ultrafilters.

\begin{lemma}
  \label{lem:Quotient ultrafilter}Let $p : X \rightarrow Y$ be a surjective
  map of sets. Then $\lambda'$ is an ultrafilter in $Y$ iff there exists an
  ultrafilter $\lambda$ in $X$ such that $p (\lambda) = \lambda'$
\end{lemma}

\begin{proof}
  Suppose that there exists an ultrafilter $\lambda$ in $X$ with $p (\lambda)
  = \lambda'$. Let $V \subset Y$. If $p^{- 1} (V) \in \lambda,$then $V \in p
  (\lambda) = \lambda'$. Conversely, if $X - p^{- 1} (V) \in \lambda$, then $p
  (X - p^{- 1} (V)) = Y - V \in \lambda'$, so $\lambda'$ is an ultrafilter.
  
  Now suppose that $\lambda'$ is an ultrafilter, and consider
  \[ p^{- 1} (\lambda') \assign \{ p^{- 1} (V) \mid V \in \lambda \} \]
  Since $\lambda'$ has the finite intersection property, so does $p^{- 1}
  (\lambda')$, and therefore there is an ultrafilter $\lambda$ on $X$ which
  contains $p^{- 1} (\lambda')$. We show that $p (\lambda) = \lambda'$. First,
  since $p^{- 1} (\lambda') \subset \lambda$, it follows that $\lambda'
  \subset p (\lambda)$. Now suppose that $U \in \lambda$. If $p (U) \notin
  \lambda'$, then $Y - p (U) \in \lambda'$, and therefore $p^{- 1} (Y - p (U))
  \in p^{- 1} (\lambda') \subset \lambda$. However, $p^{- 1} (Y - p (U)) \cap
  U \subset p^{- 1} (Y - p (U)) \cap p^{- 1} (p (U)) = \emptyset$, which is a
  contradiction. therefore $p (U) \in \lambda'$, as desired.

\end{proof}

\begin{theorem}
  Let $(S^n, \tau)$ be a topological $n$-sphere with base point $\ast \in X$
  and pseudotopology $\tau$ (induced by the topological structure). Let
  $\Sigma   S^n \assign S^n \times [0, 1]  / (x, k) \sim (x', k) \sim (\ast,
  t)$, $k \in \{ 0, 1 \}, t \in [0, 1]$, denote the pseudotopological based
  suspension of $S^n$. Then, for every $n \geq 1$, $\Sigma  S^{n - 1} \cong
  \tau (\Sigma  S^{n - 1}) \cong (S^n, \tau)$, where $\tau (\Sigma S^{n - 1})$
  is the topological modification of $\Sigma S^{n - 1}$.
\end{theorem}

\begin{proof}
  We consider the quotient map $p : S^{n - 1} \times [0, 1] \rightarrow \Sigma
  S^{n - 1}$, where $S^{n - 1}$ and $[0, 1]$ are endowed with the usual
  topological pseudotopologies, and the quotient is taken in
  $\catname{PsTop}$. Note, however, that $p$ is a topological homeomorphism on
  $S^{n - 1} \times [0, 1] - Q \assign \{ (x, k) \mid k \in \{ 0, 1 \} \} \cup
  \{ (\ast, t) \mid t \in [0, 1] \}$, and that this induces a
  pseudotopological homeomorphism on the region as well, by the definition of
  the quotient structure on $\Sigma S^{n - 1}$.
  
  Now consider $[\ast, 0] \in \Sigma S^{n - 1}$, the equivalence class of the
  point $(\ast, 0) \in S^{n - 1} \times [0, 1]$, and let $\mathcal{V}$ be the
  filter of open neighborhoods of the point $[\ast, 0] \in \Sigma S^{n - 1}$
  in the standard quotient topology. We show that a filter $\lambda'$ on
  $\Sigma S^{n - 1}$ converges in the pseudotopological quotient to $[\ast,
  0]$ iff it contains  V . First, suppose that a filter $\lambda'$ in $\Sigma
  S^{n - 1}$ does not contain $\mathcal{V}$ but converges to $[\ast, 0]$. Then
  there exists an open set $V \subset \Sigma S^{n - 1}$ which is not in
  $\lambda'$. However, by definition of the quotient pseudotopology, there
  exists a filter $\lambda$ in $S^{n - 1} \times [0, 1]$ which converges to a
  point in $p^{- 1} [\ast, 0]$ such that $p (\lambda) = \lambda'$. By
  definition of the topological pseudotopology, however, $\lambda$ must
  contain $p^{- 1} (V)$, and therefore $\lambda'$ contains $V$, a
  contradiction. Therefore every filter which converges to $[\ast, 0]$ must
  contain $\mathcal{V}$.
  
  \ We now show that the filter $[\mathcal{V}]$ generated by $\mathcal{V}$ on
  $\Sigma S^{n - 1}$ converges to $[\ast, 0]$, from which it will follow that
  any filter which contains $\mathcal{V}$ converges to $[\ast, 0]$. Let
  $\lambda'$ be an ultrafilter on $\Sigma S^{n - 1}$ which contains
  $[\mathcal{V}]$. By Lemma \ref{lem:Quotient ultrafilter}, there exists an
  ultrafilter $\lambda$ on $S^{n - 1} \times [0, 1]$ such that $p (\lambda) =
  \lambda'$. Since $S^{n - 1} \times [0, 1]$ is compact, $\lambda$ converges
  to a point $(x, t) \in S^{n - 1} \times [0, 1]$.
  Suppose that $p (x, t) \neq [\ast, 0]$. Then
  $(x, t) \notin S^{n - 1} \times \{ 0, 1 \} \cup \{ \ast \} \times [0, 1]$.
  However, there are open neighborhoods $U$ of $(x, t)$ and $V$ of $S^{n - 1}
  \times \{ 0, 1 \} \cup \{ \ast \} \times [0, 1]$ such that $U \cap V =
  \emptyset$, but both $U, V \in \lambda$ by construction, which is a
  contradiction. Therefore $(x, t) \in p^{- 1} [\ast, 0]$ and $\lambda' = p
  (\lambda) \rightarrow p (x, t) = [\ast, 0]$. Since $\lambda'$ is an
  arbitrary ultrafilter which contains $[\mathcal{V}]$, it follows that
  $[\mathcal{V}] \rightarrow [\ast, 0] .$
  
  In conclusion, we have shown that the quotient pseudotopology on $\Sigma
  S^{n - 1}$ is the pseudotopology induced by the quotient topology on $\Sigma
  S^{n - 1}$, and therefore $\Sigma S^{n - 1} \cong S^n$, as desired. 
\end{proof}

\section{The Quillen Model Category Structure on
$\catname{PsTop}$}\label{sec:Model category}

In order to construct the Quillen model category structure on
$\catname{PsTop}$, we generalize the construction in {\cite{Hirschhorn_2019}}
from $\catname{Top}$ to $\catname{PsTop}$. In order to do so, we only need to
show that a critical compactness result ({\cite{Hirschhorn_2019}} Propositions
4.10) is also true for relative cell complexes in pseudotopological spaces.
Once this is accomplished, the rest of the proof in {\cite{Hirschhorn_2019}}
then applies verbatim to $\catname{PsTop}$, with the change that all of the
maps, diagrams, and operations take place in $\catname{PsTop}$ and not in
$\catname{Top}$. In the following, we begin with a discussion of compactness
in $\catname{PsTop}$, after which we proceed to show the necessary
generalization of {\cite{Hirschhorn_2019}}, Proposition 4.10, which completes
the construction, following {\cite{Hirschhorn_2019}}, of the Quillen model
category in the pseudotopological case. We do not repeat the full argument of
{\cite{Hirschhorn_2019}}
here.

\subsection{Compactness in PsTop}

We state the definition and several critical results for compactness in
$\catname{PsTop}$.

\begin{definition}
  A pseudotopological space $(X, \Lambda)$ is compact iff every ultrafilter on
  $X$ converges.
\end{definition}

Our first observation is that the image of a compact convergence space by a
continuous function is a compact subspace.

\begin{proposition}[{\cite{Beattie_Butzmann_2002}}, Proposition 1.4.7]
  Let $(X, \Lambda_X)$ be a compact pseudotopological space and let $(Y,
  \Lambda_Y)$ be any pseudotopological space. If $f : (X, \Lambda_X)
  \rightarrow (Y, \Lambda_Y)$ is a continuous surjection, then $(Y,
  \Lambda_Y)$ is compact.
\end{proposition}

As in topological and closure spaces, there is a characterization of
compactness in terms of coverings. The relevant notion of covering in PsTop is
the following.

\begin{definition}
  Let $(X, \mathcal{F})$ be a pseudotopological space, and suppose that $A
  \subset X$. We say that a collection $\mathcal{C}$ of sets is a
  \emph{covering system} of $A$ iff for every filter $\lambda \rightarrow
  x \in A$, there exists a set $C \in \mathcal{C}$ such that $C \in \lambda$.
  If $A = \{ x \}$, then we say that $\mathcal{C}$ is a local covering system
  of $(X, \mathcal{F})$ at the point $x$.
\end{definition}

It is important to note, however, that covering systems are distinct from
covers of sets. We recall the definition of a cover to emphasize this point.

\begin{definition}
  We say that a collection of subsets $\mathcal{C}$ of a set $X$ is a cover of
  $X$ iff $\cup_{U \in \mathcal{C}} U = X$. If  C  is a cover $\mathcal{C}$ of
  $X$ and a subcollection $\mathcal{C}' \subset \mathcal{C}$ is also covers
  $X$, then we call $\mathcal{C}'$ a \emph{subcover} of  C .
\end{definition}

We now define the interior of a collection of subsets of a pseudotopological
space.

\begin{definition}
  Let $(X, \Lambda)$ be a pseudotopological space, and let $\mathcal{U}$ be a
  collection of subsets of $X$. We define the interior of $\mathcal{U}$,
  $\tmop{int} (\mathcal{U})$, to be the set
  \[ \tmop{int} (\mathcal{U}) \assign \left\{ x \in X \mid \mathcal{U}  \text{
     is a covering system of } (X, \Lambda)  \text{at } x \right\} \]
  When $\mathcal{U} = \{ U \}$, we will write $\tmop{int} (U)$ or $\ring{U}$
  for $\tmop{int} (\mathcal{U})$.
\end{definition}

\begin{proposition}
  If $(X, \Lambda_X)$ is a pseudotopological space and $(Y, \Lambda_Y) \subset
  (X, \Lambda_X)$ is a topological subspace, then the topological interior of
  $Y$ is equal to $\tmop{int}_{\Lambda_X} (Y)$.
\end{proposition}

We also recall the generalization to $\catname{PsTop}$ of the property that
closed subsets of a compact set are compact. We make this precise with the
following.

\begin{definition}
  Let $(X, \Lambda)$ be a pseudotopological space. For each subset $A \subset
  X$, we define the \emph{adherence} $a_{\Lambda} (A) \assign \{ x \in X
  \mid \exists (\lambda, x) \in \Lambda \tmop{such} \tmop{that} A \in \lambda
  \}$. We say that $A \subset X$ is \emph{closed} iff $a_{\Lambda} (A) =
  A$.
\end{definition}

\begin{proposition}[{\cite{Beattie_Butzmann_2002}}, Proposition 1.4.6$(i)$]
  \label{prop:Closed sub of compact}A closed subspace $C \subset X$ of a
  compact pseudotopological space $(X, \Lambda)$ is compact.
\end{proposition}

Finally, we recall that a compact topological space remains compact in PsTop.

\begin{proposition}
  Let $(X, \tau)$ be a compact topological space. Then $(X, \Lambda_{\tau})$
  is a compact pseudotopological space.
\end{proposition}

\begin{proposition}[{\cite{Beattie_Butzmann_2002}}, Proposition 1.4.15]
  \label{prop:Finite subcover}A pseudotopological space $(X, \mathcal{F})$ is
  compact iff every covering system $\mathcal{C}$ of $(X, \mathcal{F})$
  contains a finite subcover.
\end{proposition}

The next proposition will be needed for the proof of Theorem
\ref{thm:Compactness}.

\begin{proposition}
  \label{prop:Lcs is hereditary}Let $(X, \Lambda)$ be a pseudotopological
  space, let $U \subset X$ be a subspace of $X$, and suppose that
  $\mathcal{C}$ is a local covering system of $X$ at a point $x \in U$. Then
  the set
  \[ \mathcal{C}' \assign \{ C \cap U \mid C \in \mathcal{C} \} \]
  is a local covering system of the subspace $(U, \Lambda_U)$ at $x$.
\end{proposition}

\begin{proof}
  Recall that $ (\lambda, x) \in \Lambda_U$ iff $([\lambda], x) \in \Lambda$,
  where $[\lambda]$ is the filter in $X$ generated by the sets in $\lambda$.
  Now suppose that  C  is a local covering system of $(X, \Lambda)$ at a point
  $x \in U$, and define $\mathcal{C}'$ as in the statement of the proposition.
  Then for every $(\lambda, x) \in \Lambda_U$, $\mathcal{C}$ contains a set $A
  \in [\lambda]$. Since $[\lambda]$ is generated by $\lambda$, there exists $B
  \in \lambda$ such that $B \subset A$. However, $B \subset U$ by
  construction, and therefore $B \subset A \cap U.$ It follows that $A \cap U
  \in \lambda \cap \mathcal{C}'$. Since $(\lambda, x)$ is an arbitrary element
  of $\Lambda_U$ which converges to $x \in U$, $\mathcal{C}'$ is a local
  covering system of $(U, \Lambda_U)$ at $x \in U$, and the proof is complete.
\end{proof}

We now state and prove the compactness theorem necessary for the construction
of the Quillen model category structure on $\catname{PsTop}$. We begin with
the definition of a cell attachment and a relative cell complex.

\begin{definition}
  \label{def:Cell attachment}If $X$ is a subspace of $Y$ and there is a
  pushout square
  \[
     \begin{tikzcd}
     S^{n-1} \ar[r] \ar[d] & X
     \ar[d,"p"]\\
     D^n\ar[ru,dashed,"h"]\ar[r] & Y
    \end{tikzcd}
\]
  for some $n \geq 0$, then we say that $Y$ is obtained from $X$ by
  \emph{attaching a cell} or a \emph{cell attachment}. When $n = 0$,
  we let $S^{- 1} = \emptyset$ and $D^0 = \{ \ast \}$, the one-point set.
\end{definition}

\begin{definition}
  We say that a continuous map $f : X \rightarrow Y$ between pseudotopological
  spaces is a \emph{relative cell complex} if $f$ is an inclusion and $Y$
  can be constructed from $X$ by a (possibly infinite, and even transfinite)
  sequence of cell attachments as in Definition \ref{def:Cell attachment}. If
  $Y$ may be constructed from $X$ by attaching a finite number of cells, then
  we say that it is a \emph{finite relative cell complex}, and if
  $\emptyset \rightarrow Y$ is a \emph{(finite) relative cell complex},
  then we say that $Y$ is a \emph{(finite) cell complex}.
\end{definition}

\begin{definition}
  Given a relative cell complex $f : X \rightarrow Y$, a (possibly
  transfinite) sequence $(e_0 = f, e_1, \ldots )$ of cell attachments which
  constructs $Y$ from $X$ is called a \emph{presentation} of $f$. Given a
  presentation $\mathcal{E} \assign (e_0 = f, e_1, \ldots)$ of a relative cell
  complex $f : X \rightarrow Y$, we say that an ordinal number $\gamma$ is the
  \emph{presentation ordinal} of $e$ in $\mathcal{E}$ iff e is the
  $\gamma$-th cell attachment in $\mathcal{E}$. In particular, the
  presentation ordinal of $X$ is $0$. We say that $\gamma$ is the presentation
  ordinal of $f : X \rightarrow Y$ iff $\mathcal{E} \cong \gamma$.
\end{definition}

\begin{remark}
  Note that any two presentations $\mathcal{E}$ and $\mathcal{E}'$ of $f : X
  \rightarrow Y$ have the same presentation ordinal. Otherwise, without loss
  of generality, there is a cell in $\mathcal{E}$ which is not in
  $\mathcal{E}'$, so they cannot be constructions of the same space $Y$ from
  $X$.
\end{remark}

We now state and prove the generalization of {\cite{Hirschhorn_2019}},
Proposition 4.10 to $\catname{PsTop}$.

\begin{theorem}
  \label{thm:Compactness}If $X \rightarrow Y$ is a relative cell complex in
  PsTop, then a compact subset of $Y$ intersects the interiors of only
  finitely many cells of $Y - X$.
\end{theorem}

\begin{proof}
  We follow the general strategy of the proof of {\cite{Hirschhorn_2019}},
  Proposition 4.10.
  
  Let $C$ be a compact subset of $Y$. We construct a subset $P$ of $C$ by
  choosing a point of $C \cap \ring{D}$ for each cell $D$ such that the
  intersection $C \cap \ring{D}$ is non-empty. We will first show that $P$ is
  closed in $C$, and therefore compact by Proposition \ref{prop:Closed sub of
  compact}. Suppose that $c \in C$. The first step to prove that $P$ is closed
  in $C$ will be to show that there is a local covering system $\mathcal{U}_c$
  of $Y$ at $c$ such that
  \[ \bigcup_{U \in \mathcal{U}} (U \cap P) = \left\{\begin{array}{ll}
       \emptyset & c \notin P,\\
       c & c \in P.
     \end{array}\right. \]

  If $c \in C \cap (Y - X)$, then let $e_c$ be the unique cell of $Y - X$
  which contains $c$ in its interior. Otherwise, if $c \in C \cap X$, let $e_c
  = X$. Since there is at most one point of $P$ in the interior of any cell of
  $Y - X$ and the cells are homeomorphic to the topological space $D^n$, we
  may choose an local $(e_c, \Lambda_{e_c})$-covering system (i.e. a local
  covering system in the subspace $(e_c, \Lambda_{e_c}) \subset (Y,
  \Lambda_Y)$) $\mathcal{V }_c = \{ V_c \}$ of $c$ consisting of a single set
  subset of $Y$ such that $V_c \subset \tmop{int} (e_c)$ and which satisfies
  \[ V_c \cap P = \left\{\begin{array}{ll}
       c & c \in P,\\
       \emptyset & c \notin P.
     \end{array}\right. \]
  In particular, if $e_c$ is a cell of $Y - X$, then we may choose $V_c$ to be
  a small open neighborhood of $c$ which avoids any point $p \in  \ring{e}_c
  \cap P$ where $p \neq c$. If $e_c = X$, then we may take $V_c = X$. Since $P
  \cap X = \emptyset$ and $\ring{X} = X$ by definition, the condition is
  satisfied.
  
  While $c$ is only in the interior of a single cell $e_c$, $c$ may be glued
  to the boundary of subsequent cells in the process of constructing the
  relative cell complex. We will now use Zorn's lemma to show that we can
  extend $\mathcal{V}_c$ to a local $(Y, \Lambda_Y) $-covering system of $c$,
  which we call $\mathcal{U}_c$, such that
  \[ \bigcup_{U \in \mathcal{U}_c} U \cap P = \left\{\begin{array}{ll}
       c & c \in P,\\
       \emptyset & c \notin P.
     \end{array}\right. \]
  Fix a presentation of $X \rightarrow Y$, and let $\alpha$ be the
  presentation ordinal of the cell $e_c$ in the presentation, where, in
  particular, $\alpha = 0$ if $e_c = X$. Suppose that the presentation ordinal
  of the relative cell complex $X \rightarrow Y$ is $\gamma$. Consider the set
  $T$ of ordered pairs $(\beta, \mathcal{W})$, where $\alpha \leq \beta \leq
  \gamma$ and $\mathcal{W}$ is a local covering system of $c$ in $Y^{\beta}$
  such that
  \[ \mathcal{W} \cap Y^{\alpha} \assign \{ Y^{\alpha} \cap W \mid W \in
     \mathcal{W} \} = \mathcal{V}_c, \]
  where $Y^{\alpha}$ is the result of attaching the first $\alpha$ cells to
  $X$ in the presentation of $X \rightarrow Y$. We define a preorder on $T$ by
  $(\beta_1, \mathcal{W}_1) < (\beta_2, \mathcal{W}_2)$ iff $\beta_1 <
  \beta_2$ and $\mathcal{W}_2 \cap Y^{\beta_1} = \mathcal{W}_1$. If $\{
  (\beta_s, \mathcal{W}_s) \}_{s \in S}$ is a chain in $T$, then $(\cup_{s \in
  S} \beta_s, \cup_{s \in S} \mathcal{W}_s)$ is an upper bound in $T$ of the
  chain, so by Zorn's lemma, we have that $T$ has a maximal element, which we
  denote $(\alpha, \mathcal{A})$.
  
  We now show that $\alpha = \gamma$. Suppose instead that $\alpha < \gamma$,
  and consider the cells of presentation ordinal $\alpha + 1$. Since $Y$ has
  the final pseudotopology determined by $X$ and the cells of $Y - X$, it is
  enough to enlarge each $A \in \mathcal{A}$ to an $A' \subset Y^{\alpha + 1}$
  so that the intersection of $A'$ with the cell of presentation ordinal
  $\alpha + 1$ is open in that cell, and so that $A \cap P = A' \cap P$
  (which, in particular is equal to either $\emptyset$ or $c$). If $h : S^{n -
  1} \rightarrow Y^{\alpha + 1}$ is the attaching map for the cell $e$ of
  presentation ordinal $\alpha + 1$, then $h^{- 1} (A)$ contains an open set
  $B_{A, e} \subset S^{n - 1}$ for at least one $A \in \mathcal{A}$. Let
  $B_e'$ be an open collar neighborhood of $B$, which, in the case that $e
  \cap P \neq \emptyset$, is chosen to avoid the unique point $p \in e \cap
  P$. For each $A \in \mathcal{A},$ let $A'$ be the union of $A$ with all of
  the $B_{A, e}'$. Then define $\mathcal{A}' \assign \{ A' \mid A \in
  \mathcal{A} \}$. It follows that the pair $(\alpha + 1, \mathcal{A}')$ is an
  element of $T$ greater than $(\alpha, \mathcal{A})$, contradicting the
  maximality of $(\alpha, \mathcal{A})$. We therefore have that $\alpha =
  \gamma$, and we define $\mathcal{U}_c \assign \mathcal{A}$.
  
  Now consider the local covering system $\mathcal{U} \assign \cup_{c \in C}
  \mathcal{U}_c$ of $C$ in $Y$. If a point $c \in C, c \notin P$ is in the
  adherence $a_C (P)$, then there exists a filter $\lambda \rightarrow c$ in
  $C \subset Y$ such that $P \in \lambda$, however, this contradicts the
  construction of $\mathcal{U}_c$, since $P \cap U = \emptyset$ for all $U \in
  \mathcal{U}_c .$ Therefore, $a_C (P) = P$, and $P$ is closed in $C \subset
  Y$. $P$ is therefore compact in $C$. We also see that, when $c \in P$, the
  above procedure produces a $\Lambda_Y$-local covering system $\mathcal{U}_c$
  of $c$ in $Y$ such that $\cup_{U \in \mathcal{U}_c} (U \cap P) = \{ c \}$.
  Taking the union $\mathcal{U} : = \cup_{c \in C} \mathcal{U}_c$, we have a
  $\Lambda_Y$-covering system of $C$, each of whose sets contains at most one
  point of $P$. By Proposition \ref{prop:Finite subcover}, there exists a
  finite subcover $\mathcal{U}' \subset \mathcal{U}$ of $C$, and since $P
  \subset C$, we conclude that $P$ has at most a finite number of points. By
  the construction of $P$, it then follows that $C$ intersects the interiors
  of at most a finite number of cells of $Y - X$.
\end{proof}

\subsection{The Quillen Model Category Structure}

We begin this section by recalling the definition of a model category
structure and state the main theorem of this section. As in
{\cite{Hirschhorn_2019}}, we use the definition of a model category from
{\cite{Hirschhorn_2003}}, Definition 7.1.3, in which a model category is
required to admit all small (co)limits instead of only finite ones, and we ask
that the factorizations not only exist but are also functorial. We first
recall the definition of a model category from {\cite{Hirschhorn_2003}}.

\begin{definition}
  A \emph{model category} is a category  M  with three distinguished
  classes of morphisms, called \emph{weak equivalences, cofibrations,} and
  \emph{vibrations}, which satisfy the following axioms:
  \begin{enumerate}
    \item (Limits) The category  M  is complete and co complete.
    
    \item (2-out-of-3) If $f$ and $g$ are maps in  M  such that $g f$ is
    defined and any two of $f, g,$ and $g f$ are weak equivalences, then so is
    the third.
    
    \item (Retract) If $f$ and $g$ are maps in  M  such that $f$ is a retract
    of $g$ and $g$ is a weak equivalence, a vibrations, or a cofibration, then
    so is $f$.
    
    \item (Lifting) Given the commutative solid arrow diagram
   	\[
       \begin{tikzcd}
       A \ar[d,"i"]\ar[r] & X \ar[d,"p"]\\
       B \ar[r] \ar[ur,dashed] & Y
       \end{tikzcd}
   \]
    the dashed arrow exists if either
    \begin{enumerate}
      \item $i$ is a cofibration and $p$ is a trivial fibration, or
      
      \item $i$ is a trivial fibration and $p$ is a fibration.
    \end{enumerate}
    \item (Factorization) Every morphism $g$ in  M  has two functorial
    factorizations:
    \begin{enumerate}
      \item $g = \rho (g) \circ \gamma (g)$, where $\rho (g)$ is a trivial
      fibration and $\gamma (g)$ is a cofibration, and
      
      \item $g = \beta (g) \circ \iota (g)$, where $\iota (g)$ is a trivial
      cofibration and $\beta (g)$ is a fibration. \ 
    \end{enumerate}
  \end{enumerate}
\end{definition}

We now define the notions on pseudotopological spaces with which we will
construct the Quillen model category on $\catname{PsTop}$. We will typically
refer to the pseudotopological space $(X, \Lambda)$ as $X$ when the structure
$\Lambda$ is unambiguous.

\begin{definition}
  A continuous function $f : X \rightarrow Y$ between pseudotopological spaces
  is a \emph{weak homotopy equivalence} iff either $X$ and $Y$ are both
  empty or $X$ and $Y$ are both non-empty and the induced map
  \[ f_{\ast} : \pi_n (X, x) \rightarrow \pi_n (Y, f (x)) \]
  is an isomorphism for every $n \geq 0$ and every choice of base-point $x \in
  X$.
\end{definition}

\begin{definition}
  A map $\alpha$ is said to be a \emph{retract} of a map $\beta$ iff there
  exists a commutative diagram
     \begin{equation*}
     \begin{tikzcd}
     A \ar[rr,bend left,"1_A"] \ar[d,"\alpha"] \ar[r] & C \ar[r] \ar[d,"\beta"] & A \ar[d,"\alpha"]\\
     B \ar[r]\ar[rr,bend right,"1_B"'] & D \ar[r] & B
    \end{tikzcd} 
     \end{equation*}
\end{definition}

\begin{definition}
  We say that a function $p : X \rightarrow Y$ is a \emph{Serre fibration}
  iff any diagram of solid arrows
  \begin{equation*}
     \begin{tikzcd} 
     D^n \ar[r] \ar[d] & X \ar[d,"p"]\\     
     D^n \times I \ar[ru,dashed,"h"]\ar[r] & Y
    \end{tikzcd}
  \end{equation*}
  admits a lift $h : D^n \times I \rightarrow X$.
\end{definition}

We now state and prove the main theorem of this section.

\begin{theorem}
  There is a model category structure on the category PsTop of
  pseudotopological spaces where
  \begin{enumerate}
    \item the weak equivalences of the model category are weak homotopy
    equivalences,
    
    \item the cofibrations of the model category are continuous functions
    which are either relative cell complexes or retracts of relative cell
    complexes, and
    
    \item the fibrations of the model category are Serre fibrations.
  \end{enumerate}
\end{theorem}

\begin{proof}
  The proof is identical to the proof of {\cite{Hirschhorn_2019}}, Theorem
  2.5, with the exception that Theorem \ref{thm:Compactness} replaces
  {\cite{Hirschhorn_2019}}, Proposition 4.10, and that all morphisms and
  constructions in the proof are in the category PsTop instead of in Top.
\end{proof}

\subsection{Fibrant and Cofibrant Objects}

In this final section, we show that all objects in $\catname{PsTop}$
are fibrant and weakly homotopy equivalent to either a cell complex or the
retract of a cell complex.

\begin{definition}
  A \emph{fibrant} object in a model category $\mathcal{C}$ with a
  terminal object $\{ \ast \}$ is a  C -object $X$ such that the map $X
  \rightarrow \{ \ast \}$ is a fibration.
\end{definition}

\begin{proposition}
  Every object in $\catname{PsTop}$ is fibrant.
\end{proposition}

\begin{proof}
  The constant map $X \rightarrow \ast$ is a Serre fibration for any
  pseudotopological space $X$. 
\end{proof}

\begin{definition}
  A cofibrant object in a model category $\mathcal{C}$ with an initial object
  $\emptyset$ is a $\mathcal{C}$-object $X$ such that the map $\emptyset
  \rightarrow X$ is a cofibration.
\end{definition}

\begin{proposition}
  The cofibrant objects in the Quillen model category on
  $\catname{PsTop}$ are the cell complexes.
\end{proposition}

\begin{proof}
  Immediate from the definitions of cofibrations and cell complexes.
\end{proof}

We end with the following surprising consequence of these constructions.

\begin{theorem}
  Every pseudotopological space $X$ is weakly homotopy equivalent to a cell
  complex or the retract of a cell complex.
\end{theorem}

\begin{proof}
  The existence of the Quillen model category structure on $\catname{PsTop}$
  implies that every continuous map $g : Z \rightarrow X$ between
  pseudotopological spaces may be factored as $g = p \circ i$, where $i$ is a
  cofibration and $p$ is a fibration and a weak equivalence. In particular,
  when $Z = \emptyset$, this implies that there exists a pseudotopological
  space $Y$ such that $i : \emptyset \rightarrow Y$ is a cofibration and $p :
  Y \rightarrow X$ is a weak equivalence. By definition, $i : \emptyset
  \rightarrow Y$ is a cofibration iff $Y$ is a cell complex or the retract of
  a cell complex.
\end{proof}

\begin{remark}
  Note that although the cells in a cell complex are topological, the process
  of construction a cell complex in $\catname{PsTop}$ may a priori
  result in a non-topological space, although we are unaware of an example
  where this occurs.
\end{remark}

\section*{Acknowledgments}

We thank Chris Kapulkin, Daniel Carranza, and Omar Antol{\'i}n Camarena for
their questions and comments in discussions on this work, as well as for their
encouragement and enthusiasm for the questions treated here. We are also
grateful to MSRI, the Unidad Cuernavaca of the Instituto de Matem{\'a}ticas de
la Universidad Nacional Aut{\'o}noma de M{\'e}xico, and the organizers of the
MSRI program ``Higher Categories and Categorification, Part Two'' for
supporting my participation in the MSRI program in Cuernavaca during June
2022, and for the extremely pleasant working atmosphere they provided there.

\end{document}